\documentclass[final,leqno]{siamltex}
\usepackage{graphicx}
\usepackage{amsmath}
\usepackage{amsfonts}
\usepackage{amssymb}
\usepackage{mathrsfs}
\usepackage{verbatim}
\usepackage{subfigure}

\usepackage[notref, notcite, color]{showkeys}
\usepackage{mathtools}
\usepackage{caption}

\usepackage[letterpaper, margin=1in]{geometry}

\newcommand{\inner}[2]{#1 \cdot #2}
\newcommand{\innerb}[2]{(#1 \cdot #2)}
\newcommand{\innerbb}[2]{\left( #1 \cdot #2 \right)}
\newcommand{\innerBb}[2]{\bigl( #1 \cdot #2 \bigr)}
\newcommand{\eps}{\epsilon}

\newcommand{\nono}{\nonumber}
\newcommand{\ter}[1]{\textcolor{red}{#1}}

\newcommand{\Er}[1]{\mathcal{E}_{#1}}

\renewcommand{\grad}{\nabla}
\newcommand{\vp}{\vec{p}}
\newcommand{\vx}{\vec{x}}
\newcommand{\vxh}{\vec{x}^h}
\newcommand{\vxm}{|\vx_\rho|}
\newcommand{\vxhm}{|\vec{x}_\rho^h|}
\newcommand{\thh}{\vec{E}}
\newcommand{\vxi}{\vec{\xi}}
\newcommand{\vtau}{\vec{\tau}}
\newcommand{\vnu}{\vec{\nu}}
\newcommand{\vtauh}{\vec{\tau}^h}
\newcommand{\vnuh}{\vec{\nu}^h}
\newcommand{\wh}{w^h}
\newcommand{\zh}{Z}
\newcommand{\vX}{\vec{X}}

\newcommand{\vTau}{\vec{\mathcal{T}}}
\newcommand{\vNu}{\vec{\mathcal{V}}}
\newcommand{\vXmn}[1]{|\vX_\rho^{#1}|}
\newcommand{\vE}{\vec{E}}
\newcommand{\egs}{e^{-\gamma s} \,}
\newcommand{\egt}{e^{-\gamma t} \,}
\newcommand{\ints}{\int_0^t}
\newcommand{\intsT}{\int_0^{\Ts}}
\newcommand{\hs}{h^\star}

\newcommand{\Ts}{T_h^\star}
\newcommand{\Th}{T_h}
\newcommand{\zi}{0,\infty}
\newcommand{\oi}{1,\infty}

\newcommand{\ds}{\; {\rm d}s}
\newcommand{\drho}{\; {\rm d}\rho}

\newcommand{\Ip}{\mathcal{I}}
\newcommand{\Iip}{(0,1)}

\newcommand{\Dt}{\Delta t}
\newcommand{\gd}{w_b}

\def\epsilon{\varepsilon}

\title{Numerical analysis for a system coupling curve evolution attached orthogonally to a fixed boundary, to a reaction-diffusion equation on the curve}
\author{Vanessa Styles\footnotemark[1] \and James Van Yperen\footnotemark[1]}

\date{~}
\begin{document}

\maketitle

\begin{abstract}
We consider a semi-discrete finite element approximation for a system
consisting of the evolution of a planar curve evolving by forced curve shortening flow inside a given bounded domain $\Omega \subset \mathbb{R}^2$, such that the curve meets the boundary $\partial\Omega$ orthogonally, and the forcing is a function of the solution of a reaction-diffusion equation that holds on the evolving curve. 
 We prove optimal error bounds for the resulting approximation and present numerical experiments. 
\end{abstract}




\footnotetext[1]{Department of Mathematics, University of Sussex, Brighton, BN1 9RF, UK}

\section{Introduction}

We consider a family of planar curves $\Gamma(t)$ evolving by forced curve shortening flow inside a given bounded domain $\Omega \subset \mathbb{R}^2$, such that the curve meets the boundary $\partial\Omega$ orthogonally, and the forcing is a function of the solution of a reaction-diffusion equation that holds on $\Gamma(t)$. We combine the parametrisation presented in \cite{EF16}, for the setting in which $\Gamma(t)$ is a closed curve, with the parametrisation presented in \cite{DE98}, for the setting in which $\Gamma(t)$ meets the boundary $\partial\Omega$ orthogonally, to yield the following system of partial differential equations:  \\
find $\vx \text{ : } [0,1]\times[0,T]\rightarrow\mathbb{R}^2$ and $w \text{ : } [0,1]\times[0,T]\rightarrow\mathbb{R}$ such that
\begin{eqnarray}
&\alpha \vx_t + (1 - \alpha) \innerb{\vx_t}{\vnu} \vnu - \frac{\vx_{\rho\rho}}{\vxm^2} = \, f(w) \, \vnu & \text{in } \Iip \times (0,T), \label{xeq} \\
&(\vxm \, w)_t - \left( \frac{w_\rho}{\vxm} \right)_\rho - (\psi \, w)_\rho = \, \vxm \, g(v,w) & \text{in } \Iip \times (0,T), \label{weq} \\
&F(\vx(0,t)) = F(\vx(1,t)) = 0 & t \in [0,T], \label{Feq}\\
&\inner{\vx_\rho(0,t)}{\nabla^\perp F(\vx(0,t))} = \inner{\vx_\rho(1,t)}{\nabla^\perp F(\vx(1,t))} = 0  \qquad  &t \in [0,T], \label{xFp} \\
&w(0,t) =   w(1,t) = \gd & t \in [0,T], \label{wbc} \\
&\vx(\rho,0) = \vx^{\, 0}(\rho),~~w(\rho,0) = w^0(\rho) & \rho \in \Iip. \label{xic}
\end{eqnarray}
Here $\alpha \in (0,1]$, $\vx(\cdot,t)$ denotes the parametrisation of $\Gamma(t)$ with $\vx^{\, 0}$ parametrising the initial curve $\Gamma(0)$, $\vtau$ and $\vnu$ respectively denote unit tangent and normal vectors of $\Gamma(t)$ such that
\begin{align}
\vtau = \frac{\vx_\rho}{\vxm}, \qquad \vnu = \vtau^\perp \label{tn}
\end{align}
where, for some $\vp \in \mathbb{R}^2$, we fix $(\vp_0,\vp_1)^\perp = (-\vp_1,\vp_0)$, $\psi$ and $v$ respectively denote the tangential and normal velocities of $\Gamma(t)$
\begin{align}
\SwapAboveDisplaySkip
\psi = \inner{\vx_t}{\vtau}, \qquad v = \inner{\vx_t}{\vnu}, \label{vp}
\end{align}
and we assume that $\partial \Omega$ is given by a smooth function $F$ such that
\begin{align*}
\partial \Omega = \{\vp \in \mathbb{R}^2 : F(\vp) = 0\}
\end{align*}
which in addition we assume satisfies
\begin{align}
\SwapAboveDisplaySkip
|\nabla F(\vp)| = 1 \qquad \vp \in \partial\Omega. \label{norm_F}
\end{align}

For a closed curve $\Gamma(t)$ in $\mathbb{R}^2$, the formulation of curve shortening flow, in the form of \eqref{xeq} with $f(w) = 0$, was presented and analysed in \cite{EF16}, where the DeTurck trick is used in coupling the motion of the curve to the harmonic map heat flow with the parameter $\alpha\in (0,1]$ being such that $1/\alpha$ corresponds to the diffusion coefficient in the harmonic map heat flow. Setting $\alpha\in (0,1]$ introduces a tangential component in the velocity which, at the numerical level, gives rise to a good distribution of the mesh points along the curve. Setting $\alpha =1$ one recovers the formulation introduced and analysed in \cite{DD95}, while formally setting $\alpha = 0$ yields the approach introduced in \cite{BGN11}. The associated closed curve formulation of \eqref{xeq}--\eqref{xic} was studied in \cite{BDS17}, in which the authors proved optimal error bounds for a fully discrete finite element approximation of the coupled system. While in \cite{PS15} an alternative formulation, again for closed curves in $\mathbb{R}^2$, was presented together with optimal error bounds for a semi-discrete finite element approximation of the coupled system. Setting $\alpha=1$ and $f(w)=0$ in \eqref{xeq} and coupling the resulting equation to \eqref{xFp}, \eqref{Feq} gives rise to the model presented and analysed in \cite{DE98}, in which optimal order error bounds for a semi-discrete finite element approximation of curve shortening flow with a prescribed normal contact to a fixed boundary are presented.

The coupled system \eqref{xeq}--\eqref{xic}, and the associated closed curve formulation studied in \cite{BDS17}, can both be used to model diffusion induced grain boundary motion, \cite{H88}. This phenomenon can be observed if a polycrystalline film of metal is placed in a vapour containing another metal: atoms from the vapour diffuse into the film along the grain boundaries that separate the crystals in the film, this gives rise to variations of elastic energy in the film that cause the grain boundaries to move. Physically $\Omega(t)$ represents the polycrystalline film, $\Gamma(t)$ represents a grain boundary and $w$ represents the concentration of atoms from the vapour. The closed curve formulation arises from the physical set-up in which the polycrystalline film is assumed to be very thin in the $x_3-$direction, and the resulting two-dimensional problem is obtained by assuming independence in the $x_3-$direction. While in the set-up we consider here, the film is assumed to be infinitely long in the $x_2-$direction such that the resulting two-dimensional problem is obtained by assuming independence in the $x_2-$direction, and the grain boundary is assumed to span the width ($x_3$-direction) of the film such that it meets the upper and lower surfaces of the film orthogonally.  A more in-depth derivation of the physical set-up can be found in \cite{H88} and \cite{MS99}.

\section{Weak formulation and finite element approximation}

\subsection{Notation for function spaces}

We set $\Ip=(0,1)$ and adopt the standard notation for Sobolev spaces $W^{l,p}(\Ip)$, where $l \in \mathbb{N}_0$ and $p \in [1,\infty]$, with the Sobolev $l,p$ norm of a function $f$ on the interval $\Ip$ defined to be $\|f\|_{l,p}$ and its associated seminorm to be $|f|_{l,p}$. For the special case of $p=2$, we denote $W^{l,2}(\Ip)$ by $H^l(\Ip)$ with the associated norm and seminorm denoted by $\|f\|_{l}$ and $|f|_{l}$ respectively.  When the function is vector-valued, the function spaces are naturally extended to $[W^{l,p}(\Ip)]^n$ and $[H^l(\Ip)]^n$ with appropriately defined norms and seminorms, we, however, leave the notation for the norms unchanged. We extend the notation to include time-dependent spaces $W^{l,p}(I;X)$, where $I \subset \mathbb{R}$ is the time-dependent space and $X$ is a Banach space, with the standard associated norm and seminorm $\|f\|_{W^{l,p}(I;X)}$ and $|f|_{W^{l,p}(I;X)}$ respectively. We denote the $L^2(\Ip)$-inner product by $(f,g)$.

\subsection{Weak formulation}

Multiplying \eqref{xeq} by $\vxi \, \vxm^2$, where $\vxi \in [H^1(\Ip)]^2$ is a test function, integrating over $\Ip$ and using integration by parts gives
\begin{align}
    \SwapAboveDisplaySkip
\left( \vxm^2 \,(\alpha \, \vx_t+(1-\alpha)(\vx_t\cdot\vnu) \, \vnu),\vxi \right) + \left(\vx_\rho,\vxi_\rho \right)= \left[ \inner{\vx_\rho}{\vxi} \right]_0^1 + \left( \vxm^2 f(w) \, \vnu,{\vxi} \right). \label{prewf}
\end{align}
Using \eqref{xFp} and \eqref{norm_F} we have
\begin{align*}
\inner{\vx_\rho}{\vxi}  &= \inner{\Bigl(\nabla F(\vx) \innerb{\vx_\rho}{\nabla F(\vx)}+ \nabla^\perp F(\vx) \innerb{ \vx_\rho}{\nabla^\perp F(\vx)} \Bigr)}{ \Bigl(\nabla F(\vx) \innerb{\vxi}{\nabla F(\vx)} + \nabla^\perp F(\vx) \innerb{\vxi}{\nabla^\perp F(\vx)}\Bigr)} \\
&= \inner{\Bigl(\nabla F(\vx) \innerb{\vx_\rho}{\nabla F(\vx)}\Bigr)}{ \Bigl(\nabla F(\vx) \innerb{\vxi}{\nabla F(\vx)} + \nabla^\perp F(\vx) \innerb{\vxi}{\nabla^\perp F(\vx)}\Bigr)} \\
&= \innerb{\grad F(\vx)}{\grad F(\vx)}\innerb{\vx_\rho}{\nabla F(\vx)} \innerb{\vxi}{\nabla F(\vx)} + \innerb{\nabla F(\vx)}{\nabla^\perp F(\vx)} \innerb{\vx_\rho}{\nabla F(\vx)} \innerb{\vxi}{ \nabla^\perp F(\vx)} \\
&= \innerb{\vx_\rho}{\nabla F(\vx)}\innerb{\vxi}{\nabla F(\vx)}
\end{align*}
which combined with (\ref{prewf}) yields the following weak formulation of \eqref{xeq},  \eqref{xFp}: for all $\vxi \in [H^1(\Ip)]^2$
\begin{align}
\left( \vxm^2 \left[\alpha \,\vx_t +(1-\alpha)(\vx_t\cdot\vnu) \, \vnu \right], \vxi \, \right) + \left( \vx_\rho, \vxi_\rho \right) = \left[ \innerb{\vx_\rho}{\nabla F(\vx)}\innerb{\vxi}{\nabla F(\vx)} \right]^1_0 + \left( \vxm^2 f(w) \, \vnu, \vxi \right).
\label{xwf2}
\end{align}
Multiplying \eqref{weq} by a test function $\eta \in H^1_0(\Ip)$, integrating over $\Ip$ and using integration by parts we have
\begin{align}
\label{wwf}
\left( \left(\vxm \, w\right)_t, \eta \right) + \left( \frac{w_\rho}{\vxm}, \eta_\rho \right) + \left( \psi \, w, \eta_\rho \right) = \left( \vxm \, g(v,w), \eta \right) \qquad \forall \, \eta \in H^1_0(\Ip).
\end{align}
We assume that there is a unique solution $(\vx,w)$ of \eqref{xwf2}, \eqref{wwf} on the time interval $[0,T]$ that satisfies the boundary condition (\ref{Feq}) and the initial data \eqref{xic}. Furthermore we assume that this unique solution, and the data, satisfy
\begin{align}
& \vx \in W^{1,\infty}(0,T; [H^2(\Ip)]^2) \cap W^{2,\infty}(0,T; [L^2(\Ip)]^2); \label{xreg} \\
& w \in C([0,T]; H^2(\Ip)) \cap W^{\oi}(0,T;H^1(\Ip)); \label{wreg} \\
& f \in C^{1,1}(\mathbb{R}); \label{freg} \\
& g \in C^{1,1}(\mathbb{R}^2); \label{greg} \\
& F \in C^{2,1}(\mathbb{R}^2); \label{Freg} \\
& m \leq \vxm \leq M \quad \mbox{in } [0,1] \times [0,T], \mbox{ for some $m,\,M \in \mathbb{R}_{>0}$}.\label{nodegen}
\end{align}
We note that from \eqref{xreg} and \eqref{wreg}, for any $t \in [0,T]$, we have
\begin{align}
\label{xh2}
\|\vx(\cdot,t)\|_2 + \|\vx_t(\cdot,t)\|_1 + \|\vtau(\cdot,t)\|_1 + \|\psi(\cdot,t)\|_1 + \|w(\cdot,t)\|_2 + \| w_t(\cdot,t) \|_1 \leq C.
\end{align}
We also note that, due to \eqref{Feq}, for any $t \in [0,T]$, we have
\begin{align}
\label{xtgF}
0 = \frac{d}{dt} F(\vx(\rho,t)) = \inner{\vx_t(\rho,t)}{\nabla F(\vx(\rho,t))} \qquad \text{ for } \rho \in \{0,1\}.
\end{align}

\subsection{Finite Element approximation}

We partition the interval $[0,1]$ such that $[0,1] = \cup_{j=1}^J \overline{\sigma_j}$, where $\sigma_j = (\rho_{j-1},\rho_j)$. We set $h := \max\limits_{j=1,\dots,J} h_j$, where $h_j = \rho_j - \rho_{j-1}$ and we assume that for $C > 0$
\begin{align}
\label{hdef}
h \leq C\, h_j, \qquad j=1,\dots,J.
\end{align}
We define the finite element spaces
\begin{gather}
V^h := \{ \chi^h \in C([0,1]) \text{ : } \chi^h_{|_{\sigma_j}} \text{ is affine, } j = 1,\dots,J \} \subset H^1(\Ip), \label{Vh} \\
V_0^h := \{ \chi^h \in V^h \text{ : } \chi^h(0) = \chi^h(1) = 0\} \label{Sh}
\end{gather}
and we define the basis functions of $V^h$ to be $\phi_i(\rho_j) = \delta_i^j$. We set $I^h \text{ : } C([0,1]) \rightarrow V^h$ to be the standard Lagrange interpolation operator defined as $(I^h \eta)(\rho_j) = \eta(\rho_j)$, for $j=0,\dots,J$, whereby we denote $I^h_j := I^h_{|_{\sigma_j}}$ to be the local interpolation operator. Considering $p \in (1,\infty]$, $k \in \{0,1\}$ and $l \in \{1,2\}$, the following standard interpolation results hold for $j=1,\dots,J$
\begin{align}
h_j^{\frac{1}{p}}\,|\eta^h|_{0,\infty,\sigma_j} +h_j\,|\eta^h|_{1,p,\sigma_j} &\leq C\,|\eta^h|_{0,p,\sigma_j} \qquad &&\forall \ \eta^h \in V^h,  \label{invh}\\
|(I-I^h_j)\eta|_{k,p,\sigma_j} &\leq C\,h^{\ell-k}_j\, |\eta|_{\ell,p,\sigma_j}
\qquad &&\forall \ \eta \in W^{\ell,p}(\sigma_j), \label{Ih}\\
|(I-I^h_j)\eta|_{\ell-1,\infty,\sigma_j} &\leq C\,h^{\frac{1}{2}}_j\, |\eta|_{\ell,\sigma_j} \qquad &&\forall \ \eta \in H^{\ell}(\sigma_j), \label{Ihinf}
\end{align}
where $|\eta|_{l,p,\sigma_j}$ is the seminorm of $W^{l,p}(\sigma_j)$. We define the discrete inner product $(\cdot, \cdot)^h$ and its induced norm $\| \cdot \|_h$ by
\begin{align}
\SwapAboveDisplaySkip
\left( \eta_1, \eta_2 \right)^h := \sum_{j=1}^J \int_{\sigma_j} I^h_j(\eta_1\,\eta_2) \drho, \qquad  \|\eta\|_h^2 := (\eta,\eta)^h. \label{innerh}
\end{align}
Standard interpolation theory states that, for all $\eta^h$, $\chi^h \in V^h$, and $j=1,\dots,J$, the following results hold
\begin{subequations}
\begin{align}
\int_{\sigma_j} |\eta^h|^2 \, \drho &\leq \int_{\sigma_j} I^h_j\left[|\eta^h|^2\right] \drho \leq  3\int_{\sigma_j} |\eta^h|^2 \, \drho, && \label{lumpev} \\
\left|\int_{\sigma_j} (I-I^h_j)(\eta^h\,\chi^h) \, \drho\right| &\leq C\,h_j^2 \,|\eta^h|_{1,\sigma_j}\,|\chi^h|_{1,\sigma_j} \leq C\,h_j\, |\eta^h|_{1,\sigma_j}\,|\chi^h|_{0,\sigma_j}. \label{lumperr}
\end{align}
\end{subequations}
We assign to an element $\vxh \in [V^h]^2$ a piecewise constant discrete unit tangent and normal, denoted respectively by $\vtauh$ and $\vnuh$, and on $\sigma_j$ we approximate the tangential velocity and normal velocity respectively by $\psi^h$ and $v^h$, where
\begin{align}
    \SwapAboveDisplaySkip
\vtauh = \frac{\vxh_\rho}{\vxhm}, \quad \vnuh = (\vtauh)^\perp, \quad \psi^h = \inner{\vxh_t}{\vtauh}, \quad v^h = \inner{\vxh_t}{\vnuh} \qquad ~~\mbox{on}~\sigma_j,~j=1\cdots,J. \label{tnpvh}
\end{align}
Employing the above notation we introduce the following, continuous in time, finite element approximation of (\ref{xwf2}), (\ref{wwf}): find $\vxh\text{ : } [0,1] \times [0,T] \rightarrow \mathbb{R}^2$ and $w^h \text{ : } [0,1] \times [0,T] \rightarrow \mathbb{R}$ such that $\vxh(\cdot,t) \in [V^h]^2$ and $w^h(\cdot,t) - \gd \in V_0^h$, for $t \in [0,T]$, and
\begin{align}
\biggl( \vxhm^2 \bigl[\alpha \, \vxh_t &+ (1-\alpha) \innerbb{\vxh_t}{\vnuh} \, \vnuh \bigr], \vxi^h \, \biggr)^h + \left( \vxh_\rho, \vxi^h_\rho \right) \nono \\
&= \left( \vxhm^2 f(\wh)  \,\vnuh , \vxi^h \, \right)^h + \left[ \innerbb{\vxh_\rho}{\nabla F(\vxh)} \innerBb{\vxi^h}{\nabla F(\vxh)} \right]^1_0 && \forall \, \vxi^h \in [V^h]^2, \label{xfem} \\
\left( \left( \vxhm \, \wh \right)_t ,\eta^h \right)^h &+ \left(\frac{\wh_\rho}{\vxhm},\eta^h_\rho\right) + \left(\psi^h\,\wh,\eta_\rho^h \right)^h =  \left( \vxhm\,g(v^h,\wh), \eta^h \right)^h && \forall \, \eta^h \in V_0^h, \label{wfem}
\end{align}
where $\vxh_\rho(0,t)=(\vxh_\rho)_{|_{\sigma_1}}$ and $\vxh_\rho(1,t)=(\vxh_\rho)_{|_{\sigma_J}}$, and $\vxh$, $w^h$ satisfy the boundary and initial data
\begin{align}
F(\vxh(0,t)) = F(\vxh(1,t)) = 0 \qquad t \in [0,T]&, \label{sd_bc} \\
\wh(\rho,0)=I^hw^0(\rho), ~~ \vxh(\rho,0)=I^h\vx^{\, 0}(\rho) \qquad \rho \in \Ip&. \label{sd_ic}
\end{align}
Using \eqref{Feq} and \eqref{sd_ic} we have $F(\vxh(\rho,0))=F(I^h \vx^{\,0}(\rho))=F(\vx^{\,0}(\rho))=0$, for $\rho\in\{0,1\}$, and hence the conditions \eqref{sd_bc} are equivalent to
\begin{align}
\inner{\vxh_t(0,t)}{\nabla F(\vxh(0,t))} = \inner{\vxh_t(1,t)}{\nabla F(\vxh(1,t))} = 0 \qquad & t \in [0,T], \label{sd_xtgf}
\end{align}
which is the discrete analogue of \eqref{xtgF}. Similarly we present the discrete analogue of \eqref{norm_F}, namely
\begin{align}
|\grad F(\vxh(0,t)) | = |\grad F(\vxh(1,t)) | = 1 \qquad t \in [0,T]. \label{sd_normF}
\end{align}
Let us formulate the main theorem, which will be proved in Section \ref{FEA_sec}.
\begin{theorem}
\label{thmH1fbrde}
Let $\vx^{h}(\cdot,0) = I^h \vx^{\,0}(\cdot) \in [V^h]^2$ and $w^{h}(\cdot,0) = I^h w^0(\cdot) \in V^h$. There exists $\hs > 0$ such that for all $h \in (0,\hs]$ the semi-discrete problem \eqref{xfem}-\eqref{sd_bc} has a unique solution $(\vxh,\wh) \in [V^h]^2 \times V^h$ on $[0,T]$ and the following error bounds hold
\begin{align*}
\sup_{s\in[0,T]}\left[ |\vx(\cdot,s) - \vxh(\cdot,s)|_1^2 + |w(\cdot,s) - \wh(\cdot,s)|_0^2 \right] + \int_0^T \left[ |\vx_t(\cdot,s) - \vxh_t(\cdot,s)|_0^2 + |w(\cdot,s) - \wh(\cdot,s)|_1^2 \right] \, ds \leq C h^2, \label{thmeq}
\end{align*}
for some $C>0$ independent of $h$.
\end{theorem}

\section{Error Analysis}
\label{FEA_sec}

For the proof of Theorem \ref{thmH1fbrde}, and hence throughout this section, we choose $\hs, \gamma \in \mathbb{R}_{>0}$ so that
\begin{equation}
\label{exp}
e^{\gamma \, T}(\hs)^{\frac{1}{2}} \leq \min\left\{ \frac{1}{2C_1},\beta\right\}~~\mbox{and}~~\gamma\geq \max\left\{1,\frac{32C_2}{m^2 \alpha}\right\},
\end{equation}
where $C_1,C_2\in \mathbb{R}_{>0}$ and $\beta\in(0,1]$, are independent of $h$ and will be chosen a posteriori. Standard ODE theory implies that there exists a unique solution $(\vxh, \wh)$ of \eqref{xfem}--\eqref{sd_ic} on some time interval $[0,\Th]$ $(\Th > 0)$.

For simplicity of notation we define
\begin{align*}
\thh:= I^h \vx - \vxh \quad \mbox{and} \quad \zh:= I^h w - \wh
\end{align*}
such that $\vx - \vxh = (I - I^h) \vx + \thh$ and $w - \wh= (I - I^h)w + \zh$. For the proof of Theorem \ref{thmH1fbrde} we adapt the arguments presented in \cite{DD09} and define for some $C_1 \in \mathbb{R}_{>0}$
\begin{align*}
\Ts := \sup \bigg\{ & t \in [0,T] \text{ : } (\vxh,\wh) \text{ solves  \eqref{xfem}--\eqref{sd_bc}}, \\
& \quad \frac{m}{2} \leq |\vxh_\rho| \leq 2 M \quad \text{ in } [0,1] \times [0,t], \\
& \quad \|w^h\|_{C([0,t]; L^\infty(\Ip))} \leq 2 C_w \|w\|_{C([0,T]; H^1(\Ip))}, \text{ and } \\
& \quad \sup_{s \in [0,t]} \egs \left[ |\thh(\cdot,s)|_1^2 + |\zh(\cdot,s)|_0^2 \right] + \ints \egs \left[ |\thh_t(\cdot,s)|_0^2 + |\zh(\cdot,s)|_1^2 \right] \ds < 2 C_1 h^2 \biggr\}.
\end{align*}
We then prove the result of Theorem \ref{thmH1fbrde} on $[0,\Ts]$, for $C$ independent of $\Ts$, thus enabling us to show that $\Ts = T$ and hence proving the theorem. By the definition of $\Ts$ we have the following bounds
\begin{align}
\frac{m}{2} \leq |\vxh_\rho| \leq 2 M \qquad \text{ in } [0,1] \times [0,\Ts) \label{disnodegen} \\
\|w^h\|_{C([0,\Ts); L^\infty(\Ip))} \leq 2 C_w \|w\|_{C([0,T]; H^1(\Ip))} \label{whinf} \\
\sup_{s \in [0,\Ts]} \egs \left[ |\thh(\cdot,s)|_1^2 + |\zh(\cdot,s)|_0^2 \right] + \intsT \egs \left[ |\thh_t(\cdot,s)|_0^2 + |\zh(\cdot,s)|_1^2 \right] \ds < 2 C_1 h^2. \label{cont}
\end{align}
The main part of the proof of Theorem \ref{thmH1fbrde} is split into the following two lemmas:
\begin{lemma}
\label{surf_lem}
There exists $C_2 \in \mathbb{R}_{>0}$, such that for all $t \in [0,\Ts)$, we have
\begin{align}
\label{sl}
& \frac{1}{4} \egt |\thh|_1^2 + \frac{m^2 \alpha}{16} \ints \egs |\thh_t|_0^2 \, \ds \leq C_2 \ints \egs \left[h^2 + |\thh|_1^2 + |\zh|_0^2 + h^{-1}|\thh|_{\zi}^4 \right] \ds.
\end{align}
\end{lemma}

\begin{lemma}
\label{rde_lem}
There exists $\hs>0$ and $C_3 \in \mathbb{R}_{>0}$, such that for all $h \in (0,\hs]$ and $t \in [0,\Ts)$,  we have
\begin{align}
& \frac{m}{4} \egt |\zh|_0^2  + \frac{1}{4 M} \ints \egs |\zh|_1^2 \, \ds \leq C_3 \ints \egs \left[ h^2 + |\thh_t|_0^2  + |\thh|_1^2+ |\zh|_0^2 \right] \ds.
\label{c3}
\end{align}
\end{lemma}

Before proving Lemmas \ref{surf_lem} and \ref{rde_lem} and subsequently Theorem \ref{thmH1fbrde}  we note the following useful bounds for $t \in [0,\Ts)$.

Using \eqref{Ih} and \eqref{xh2}, we have
\begin{align}
|\vx - \vxh|_1 &\leq |(I - I^h)\vx|_1 + |\thh|_1 \leq C \, h \, |\vx|_2 + |\thh|_1 \leq C\left[h + |\thh|_1 \right], \label{xxh1}
\end{align}
as well as
\begin{align}
|\vx_t - \vxh_t|_0 &\leq |(I - I^h)\vx_t|_0 + |\thh_t|_0 \leq C \, h \, |\vx_t|_1 + |\thh_t|_0 \leq C\left[h + |\thh_t|_0 \right]. \label{xtxht}
\end{align}
If we use \eqref{lumpev}, \eqref{xh2}, \eqref{invh}  and \eqref{hdef}, we get
\begin{align}
\label{xht0}
|\vxh_t|_s \leq |I^h \vx_t|_s + |\thh_t|_s \leq C \left[1 + h^{-s}|\thh_t|_0 \right] \quad \mbox{for } s=0,1.
\end{align}
Using \eqref{Ih} and \eqref{xh2}, we have
\begin{align}
|w - \wh|_0 &\leq |(I - I^h) w|_0 + |\zh|_0 \leq C \, h \, |w|_1 + |\zh|_0 \leq C \left[ h + |\zh|_0\right], \label{wwh0}
\end{align}
and from \eqref{lumpev}, \eqref{xh2}, \eqref{invh} and \eqref{hdef}, we have
\begin{align}
\label{wh1}
|\wh|_1 \leq |I^h w|_1 + C \, h^{-1} |\zh|_0 \leq C \, h^{-1} \left[ h + |\zh|_0 \right].
\end{align}
Using \eqref{nodegen}, \eqref{disnodegen} and \eqref{xxh1}, we have
\begin{align}
\label{xrxhbr}
\left| \frac{1}{\vxm} - \frac{1}{\vxhm} \right|_0 &\leq \left| \frac{\vxm - \vxhm}{\vxm \, \vxhm} \right|_0 \leq \frac{2}{m^2}|\vx - \vxh|_1 \leq C\left[h + |\thh|_1 \right].
\end{align}
In the same way, with \eqref{tn}, \eqref{tnpvh}, \eqref{nodegen} and \eqref{xxh1}, we have
\begin{align}
|\vtau - \vtauh|_0 \leq \left|\vxh_\rho \, \frac{\vxhm - \vxm}{\vxm \, \vxhm} \right|_0 + \left| \frac{1}{\vxm}(\vx_\rho - \vxh_\rho) \right|_0 \leq \frac{2}{m} \, |\vx - \vxh|_1 \leq C \left[h + |\thh|_1 \right] \nono
\end{align}
which yields
\begin{align}
\SwapAboveDisplaySkip
|\vtau - \vtauh|_0 + |\vnu - \vnuh|_0 \leq C \left[h + |\thh|_1 \right]. \label{tn0}
\end{align}
Using \eqref{vp}, \eqref{tnpvh}, \eqref{tn0}, \eqref{xtxht}, Sobolev embeddings and \eqref{xh2}, we have
\begin{align}
|\psi - \psi^h|_0 &\leq |\inner{\vx_t}{(\vtau - \vtauh)} |_0 + | \inner{\vtauh}{(\vx_t - \vxh_t)} |_0 \nono \\
&\leq C \, |\vx_t|_{0,\infty}\left[h + |\thh|_1 \right] + C \, \left[h + |\thh_t|_0 \right] \leq C \left[h + |\thh_t|_0 + |\thh|_1 \right] \nono
\end{align}
and thus we obtain
\begin{align}
\SwapAboveDisplaySkip
|\psi - \psi^h|_0 + |v - v^h|_0 \leq C \left[h + |\thh_t|_0 + |\thh|_1 \right]. \label{psipsih0}
\end{align}
With \eqref{xh2}, \eqref{invh}, \eqref{hdef} and \eqref{psipsih0}, we gain
\begin{align}
\label{vh1}
|v^h|_1 \leq |v|_1 + C \, h^{-1} \, |v - v^h|_0 \leq C \, h^{-1} \left[h + |\thh_t|_0 + |\thh|_1 \right].
\end{align}
{\bf Proof of Lemma \ref{surf_lem}: }
In this proof we combine techniques used in \cite{BDS17} and \cite{DE98}. Taking $\vxi = \thh_t$ in \eqref{xwf2} and $\vxi^h = \thh_t$ in \eqref{xfem}, subtracting the resulting equations and noting
\begin{align}
\label{xrihxr_int}
\left( \vx_\rho - (I^h \vx)_\rho, \vxi_\rho^h \right) = 0,
\end{align}
we obtain
\begin{align}
& \left( \vxhm^2 \left[\alpha \, \thh_t + (1 - \alpha) \innerb{\thh_t}{\vnuh} \, \vnuh \right], \thh_t \right)^h + \left( \thh_{\rho}, \thh_{\rho,t} \right)\nono \qquad \qquad \qquad \qquad  \\
& \qquad \qquad = \biggl[\left( \vxhm^2 \left[\alpha \, I^h\vx_t + (1-\alpha) \innerb{ I^h\vx_t}{ \vnuh} \, \vnuh \right], \thh_t \right)^h - \left( \vxm^2 \left[\alpha \,\vx_t + (1-\alpha) \innerb{\vx_t}{\vnu} \, \vnu \right], \thh_t \right) \biggr] \nono \\
& \qquad \qquad \qquad + \left[ \left( \vxm^2 f(w)  \,\vnu , \thh_t \, \right) - \left( \vxhm^2 f(\wh)  \,\vnuh , \thh_t \, \right)^h \right] \nono \\
& \qquad \qquad \qquad + \left[ \innerBb{\vx_\rho}{\nabla F(\vx)}\innerBb{\thh_t}{\nabla F(\vx)} - \innerBb{\vxh_\rho}{ \nabla F(\vxh)} \innerBb{\thh_t}{\nabla F(\vxh)} \right]^1_0 =: \sum_{i=1}^3 I_i. \label{xRHS}
\end{align}
Using \eqref{disnodegen} and \eqref{lumpev}, we note that the left-hand side of \eqref{xRHS} is bounded below
\begin{align}
& \left( \vxhm^2 \left[\alpha \, \thh_t + (1 - \alpha) \innerb{\thh_t}{\vnuh} \, \vnuh \right], \thh_t \right)^h + \left( \thh_{\rho}, \thh_{\rho,t} \right)   \nono \\
& \qquad \qquad \qquad  \geq \frac{m^2}{4} \left[ \alpha \, \|\thh_t\|_h^2 + (1 - \alpha)\|\inner{\thh_t}{\vnuh}\|_h^2 \right] + \frac{1}{2} \frac{d}{dt} |\thh|_1^2 \geq \frac{m^2 \alpha}{4} \, |\thh_t|_0^2 + \frac{1}{2} \frac{d}{dt} |\thh|_1^2. \label{xLHSb}
\end{align}
We now proceed to bound $I_1$, $I_2$ and $I_3$ in \eqref{xRHS}, beginning with $I_1$.
\begin{align}
I_1 &= \left( \vxhm^2 \left[\alpha \,  I^h\vx_t + (1-\alpha) \innerb{I^h \vx_t}{\vnuh} \, \vnuh \right], \thh_t \right)^h  - \left( \vxm^2 \left[\alpha \,\vx_t + (1-\alpha) \innerb{\vx_t}{\vnu} \, \vnu \right], \thh_t \right) \nono \\
&= \biggl[ \left( \left[\vxhm^2 - \vxm^2 \right] \left[\alpha \,\vx_t + (1-\alpha) \innerb{\vx_t}{\vnu} \, \vnu \right], \thh_t \right) \nono \\
&\qquad \quad + \, (1 - \alpha) \left( \vxhm^2 \left[ \innerb{\vx_t}{(\vnuh - \vnu)} \, \vnu + \innerb{\vx_t}{\vnuh} \left(\vnuh - \vnu \right) \right], \thh_t \right) \biggr] \nono \\
& \quad + \biggl[ \left( \vxhm^2 \left[\alpha \,  I^h\vx_t + (1-\alpha) \innerb{I^h \vx_t}{ \vnuh} \, \vnuh \right], \thh_t \, \right)^h - \left( \vxhm^2 \left[\alpha \,  I^h\vx_t + (1-\alpha) \innerb{I^h \vx_t}{\vnuh} \, \vnuh \right], \thh_t \, \right) \nono \\
& \quad + \left( \vxhm^2 \left[\alpha \,(I^h - I) \vx_t + (1-\alpha) \innerb{ (I^h - I) \vx_t}{\vnuh} \, \vnuh \right], \thh_t \, \right) \biggr] =: I_{1,1} + I_{1,2}. \label{T1}
\end{align}
Using \eqref{nodegen}, \eqref{disnodegen},  \eqref{xxh1}, \eqref{tn0},  Sobolev embeddings and \eqref{xh2}, we see that 
\begin{align}
I_{1,1} &= \left( \left[\vxhm^2 - \vxm^2 \right] \left[\alpha \,\vx_t + (1-\alpha) \innerb{\vx_t}{\vnu} \, \vnu \right], \thh_t \right) \nono \\
& \qquad + \, (1 - \alpha)  \left( \vxhm^2 \left[ \innerb{\vx_t}{(\vnuh - \vnu)} \, \vnu + \innerb{\vx_t}{\vnuh} \left(\vnuh - \vnu \right) \right], \thh_t \right) \nono \\
&\leq |\vx_t|_{\zi} \left[ \left| \vxm + \vxhm \right|_{0,\infty} | \vx - \vxh |_1 + 8 M^2 (1 - \alpha) \, |\vnu - \vnuh|_0 \right] |\thh_t|_0 \leq C \left[h + |\thh|_1  \right]|\thh_t|_0. \label{T11}
\end{align}
From (\ref{lumpev},b),  \eqref{disnodegen}, \eqref{Ih}  and \eqref{xh2}, we get %
\begin{align}
I_{1,2} &= \left( \vxhm^2 \left[\alpha \, I^h \vx_t + (1-\alpha) \innerb{I^h \vx_t}{\vnuh} \, \vnuh \right], \thh_t \, \right)^h - \left( \vxhm^2 \left[\alpha \, I^h \vx_t + (1-\alpha) \innerb{I^h \vx_t}{\vnuh} \, \vnuh \right], \thh_t \, \right) \nono \\
&\quad+\left( \vxhm^2 \left[\alpha \, (I^h - I) \vx_t + (1-\alpha) \innerb{(I^h - I) \vx_t}{\vnuh} \, \vnuh \right], \thh_t \, \right) \nono \\
&\leq C \, h  \sum_{j=1}^J \left|I^h_j \vx_t \right|_{1,\sigma_j} \left| \vxhm^2 \left[ \alpha \, \thh_t + (1 - \alpha) \innerb{\thh_t}{\vnuh} \, \vnuh \right] \right|_{0,\sigma_j}+ C \, |(I - I^h) \vx_t|_0 \, |\thh_t|_0 \nono \\
&\leq C \, h \, |\vx_t|_1 \, |\thh_t|_0 \leq C \, h \, |\thh_t|_0.
\label{T13}
\end{align}
Combining \eqref{T1}--\eqref{T13} we have
\begin{align}
|I_1| \leq \frac{m^2 \alpha}{24} \, |\thh_t|_0^2 + C \left[ h^2 + |\thh|_1^2 \right]. \label{T1b}
\end{align}
We now bound $I_2$.
\begin{align}
I_2 &= \left( \vxm^2 f(w)  \,\vnu , \thh_t \, \right) - \left( \vxhm^2 f(\wh)  \,\vnuh , \thh_t \, \right)^h \nono \\
&= \left[\left(\left[\vxm^2 - \vxhm^2 \right] f(w) \, \vnu + \vxhm^2  \, f(w) \left[ \vnu - \vnuh \right], \thh_t \right) + \left( \vxhm^2 \left[ f(w) - f(\wh) \right] \vnuh, \thh_t \right)\right] \nono \\
&\quad + \left[ \left( \vxhm^2 \left[ (I - I^h) f(\wh) \right] \vnuh, \thh_t \right) + \left( \vxhm^2 \left[ (I^h - I) f(\wh) \right] \vnuh, \thh_t \right)^h \right. \nono \\
&\quad + \left. \left( \vxhm^2 \, I^h(f(\wh)) \,\vnuh , \thh_t \, \right) - \left( \vxhm^2 \, I^h(f(\wh))  \,\vnuh , \thh_t \, \right)^h \right] =: I_{2,1} +I_{2,2}. \label{T2}
\end{align}
Using  \eqref{nodegen}, \eqref{disnodegen}, \eqref{xxh1}, \eqref{tn0}, \eqref{freg} and \eqref{wwh0}, we see that
\begin{align}
I_{2,1}&= \left(\left[\vxm^2 - \vxhm^2 \right] f(w) \, \vnu + \vxhm^2  \, f(w) \left[ \vnu - \vnuh \right], \thh_t \right)+ \left( \vxhm^2 \left[ f(w) - f(\wh) \right] \vnuh, \thh_t \right)  \nono \\
&\leq |f|_{0,\infty} \left[ \left| \vxm + \vxhm \right|_{0,\infty} |\vx_\rho - \vxh_\rho|_0 + 4 M^2 |\vnu - \vnuh|_0 \right] |\thh_t|_0 + C \, |w - \wh|_0 \, |\thh_t|_0 \nono\\
& \leq C \left[h + |\zh|_0 + |\thh|_1 \right] |\thh_t|_0. \label{T21}
\end{align}
With \eqref{disnodegen}, (\ref{lumpev},b), \eqref{Ih}, \eqref{freg}, \eqref{whinf} and \eqref{wh1}, we obtain
\begin{align}
I_{2,2}  &= \left( \vxhm^2 \left[ (I - I^h) f(\wh) \right] \vnuh, \thh_t \right) + \left( \vxhm^2 \left[ (I^h - I) f(\wh) \right] \vnuh, \thh_t \right)^h \nono \\
&\quad +  \left( \vxhm^2 \,I^h(f(\wh))  \,\vnuh , \thh_t \, \right) - \left( \vxhm^2 \,I^h(f(\wh))  \,\vnuh , \thh_t \, \right)^h \nono \\
&\leq C \, |(I - I^h)f(\wh)|_0 \,  |\thh_t|_0 +  C \, h \sum_{j=1}^J \left|I^h(f(\wh)) \right|_{1,\sigma_j} \left| \vxhm^2 \inner{\thh_t}{\vnuh} \right|_{0,\sigma_j} \nono\\
&\leq C \, h \, |f(\wh)|_1 \, |\thh_t|_0  \leq C \, h \, |f'(\wh)|_{0,\infty} \, |\wh|_1 \, |\thh_t|_0 \leq C \left[h + |\zh|_0 \right] |\thh_t|_0. \label{T23}
\end{align}
Combining \eqref{T2}--\eqref{T23}, we have
\begin{align}
|I_2| \leq \frac{m^2 \alpha}{24} \, |\thh_t|_0^2 + C \left[ h^2 + |\zh|_0^2 + |\thh|_1^2 \right]. \label{T2b}
\end{align}
We now bound $I_3$, to this end we set
\begin{align*}
b(\rho,t) := \inner{\vx_\rho}{\nabla F(\vx)}, \qquad  b^h(\rho,t) := \inner{\vxh_\rho}{\nabla F(\vxh)}
\end{align*}
and note that
\begin{align}
\label{interp_bc}
\vx(\rho,t) = I^h(\vx(\rho,t)) \quad \mbox{for } \rho \in \{0,1\}, \, t \in [0,T].
\end{align}
Using \eqref{xtgF}, \eqref{sd_xtgf}, and \eqref{interp_bc} we see that
\begin{align}
I_3 &= \left[ b(\rho,t) \innerBb{\thh_t}{\nabla F(\vx)} - b^h(\rho,t) \innerBb{\thh_t}{\nabla F(\vxh)} \right]^1_0 \nono \\
&= \bigl[ b(\rho,t) \innerBb{\thh_t}{[\nabla F(\vx) - \nabla F(\vxh)]}\bigr]_0^1 \nono \\
& \qquad \qquad  + \bigl[ (b(\rho,t) - b^h(\rho,t)) \innerBb{\vx_t}{[\nabla F(\vxh) - \nabla F(\vx)]} \bigr]^1_0  =: I_{3,1} + I_{3,2}. \label{T3}
\end{align}
Using \eqref{nodegen}, \eqref{norm_F}, Sobolev embeddings, \eqref{xreg}, \eqref{Freg}, and noting \eqref{interp_bc}, for $\rho \in \{0,1\}$ and $t \in [0,T]$ we have
\begin{subequations}
\begin{align}
|b(\rho,t)| &\leq M, \label{bbound} \\
|b_t(\rho,t)| &\leq |\vx_{\rho,t}(\rho,t)| + M | D^2 F(\vx(\rho,t)) \, \vx_t(\rho,t) | \leq C, \label{btbound}
\end{align}
as well as %
\begin{align}
\SwapAboveDisplaySkip
|\grad F(\vx(\rho,t)) - \grad F(\vxh(\rho,t))| &\leq L_{\grad F} \, |\vx(\rho,t) - \vxh(\rho,t)| \leq C \, |\thh(\cdot,t)|_{\zi}. \label{gFxgFxh}
\end{align}
\end{subequations}
A Taylor's expansion yields
\begin{align}
\nabla F(\vx) - \nabla F(\vxh) = & D^2 F(\vx) (\vx - \vxh) + \int_0^1 (D^2 F(s\vx + (1 - s)\vxh) - D^2 F(\vx))(\vx - \vxh) \, ds \label{Fexp}
\end{align}
which together with (\ref{bbound},b), \eqref{Freg}, \eqref{interp_bc}, \eqref{invh}, and \eqref{hdef}, gives
\begin{align}
I_{3,1}&= \left[ b(\rho,t) \innerBb{\thh_t}{[\nabla F(\vx) - \nabla F(\vxh)]} \right]_0^1 \nono \\
&= \left[ b(\rho,t) \, \thh^T_t D^2 F(\vx) \, \thh + b(\rho,t) \int_0^1 \thh^T_t [D^2 F(s\vx + (1 - s)\vxh) - D^2 F(\vx)]\thh \, \ds \right]_0^1 \nono \\
&= \biggl[ \frac{1}{2} \frac{d}{dt} \left( b(\rho,t) \thh^T D^2 F(\vx) \, \thh \right) - \frac{1}{2} b_t(\rho,t) \, \thh^T \, D^2 F(\vx) \, \thh - \frac{1}{2} b(\rho,t) \, \thh^T \frac{d}{dt} (D^2 F(\vx)) \, \thh \nono \\
& \qquad + b(\rho,t) \int_0^1 \thh^T_t [D^2 F(s\vx + (1 - s)\vxh) - D^2 F(\vx)]\thh \, \ds \biggr]_0^1 \nono \\
&\leq \frac{1}{2} \frac{d}{dt} \left[ b(\rho,t) \, \thh^T D^2 F(\vx) \, \thh \right]_0^1 + C \, |\thh|_{0,\infty}^2 \left[1 + |\thh_t|_{0,\infty} \right]\nono\\
&\leq \frac{1}{2} \frac{d}{dt} \left[ b(\rho,t) \, \thh^T D^2 F(\vx) \, \thh \right]_0^1 + C \, |\thh|_{0,\infty}^2 \left[1 + h^{-\frac{1}{2}}|\thh_t|_0 \right]. \label{T31}
\end{align}
Denoting $\vx(0,t) := \vx_0(t)$, taking $\vxi = (1 - \rho)\nabla F(\vx_0)$ in \eqref{xwf2} and noting \eqref{norm_F}, we have
\begin{align}
& \left( \vxm^2 \left[\alpha \,\vx_t + (1-\alpha) \innerb{\vx_t}{\vnu} \, \vnu \right], (1 - \rho)\nabla F(\vx_0) \, \right) - \left( \vx_\rho, \nabla F(\vx_0) \right) \nono \\
& \qquad = \left( \vxm^2 f(w)  \,\vnu , (1 - \rho)\nabla F(\vx_0) \, \right) + \left[ (1 - \rho) \innerBb{\vx_\rho}{ \nabla F(\vx)} \innerBb{\nabla F(\vx_0)}{\nabla F(\vx)} \right]^1_0 \nono \\
& \qquad = \left( \vxm^2 f(w)  \,\vnu , (1 - \rho)\nabla F(\vx_0) \, \right) - b(0,t), \nono
\end{align}
and hence
\begin{align}
b(0,t) &= \left( \vx_\rho, \nabla F(\vx_0) \right) + \left( \vxm^2 f(w)  \,\vnu , (1 - \rho)\nabla F(\vx_0) \, \right) \nono \\
&\quad - \left( \vxm^2 \left[\alpha \,\vx_t + (1-\alpha) \innerb{\vx_t}{\vnu} \, \vnu \right], (1 - \rho)\nabla F(\vx_0) \, \right). \label{b0}
\end{align}
Similarly, denoting $\vxh(0,t) := \vxh_0(t)$, taking $\vxi^h = (1 - \rho)\nabla F(\vxh_0)$ in \eqref{xfem} and noting \eqref{sd_normF}, we have
\begin{align}
b^h(0,t) &= \left( \vxh_\rho, \nabla F(\vxh_0) \right) + \left( \vxhm^2 f(\wh)  \,\vnuh , (1 - \rho)\nabla F(\vxh_0) \, \right)^h \nono \\
&\quad - \left( \vxhm^2 \left[\alpha \, \vxh_t + (1-\alpha) \innerb{\vxh_t}{\vnuh} \, \vnuh \right], (1 - \rho)\nabla F(\vxh_0) \, \right)^h. \label{bh0}
\end{align}
Hence, subtracting \eqref{bh0} from \eqref{b0} yields
\begin{align}
& b(0,t) - b^h(0,t) = \left[\left( \vx_\rho, \nabla F(\vx_0) \right) - \left( \vxh_\rho, \nabla F(\vxh_0) \right) \right] \nono\\
& \quad + \left[\left( \vxm^2 f(w)  \,\vnu , (1 - \rho)\nabla F(\vx_0) \, \right)  - \left( \vxhm^2 f(\wh)  \,\vnuh , (1 - \rho)\nabla F(\vxh_0) \, \right)^h \right] \nono \\
& \quad + \biggl[\left( \vxhm^2 \left[\alpha \, \vxh_t + (1-\alpha) \innerb{\vxh_t}{\vnuh} \, \vnuh \right], (1 - \rho)\nabla F(\vxh_0) \, \right)^h \nono \\
&\qquad \qquad - \left( \vxm^2 \left[\alpha \,\vx_t + (1-\alpha) \innerb{\vx_t}{\vnu} \, \vnu \right], (1 - \rho)\nabla F(\vx_0) \, \right) \biggr] =: \sum_{i=1}^3 B_i. \label{bRHS}
\end{align}
Starting with $B_1$, using \eqref{xh2}, \eqref{gFxgFxh} and \eqref{sd_normF}, and noting \eqref{xrihxr_int}, we have
\begin{align}
B_1 &= \left( \vx_\rho, \nabla F(\vx_0) \right) - \left( \vxh_\rho, \nabla F(\vxh_0) \right) \nono \\
&= \left( \vx_\rho, \nabla F(\vx_0) -  \nabla F(\vxh_0) \right) + \left( \thh_\rho, \nabla F(\vxh_0) \right) \leq C \left[ |\thh|_{0,\infty} + |\thh|_1 \right]. \label{B1}
\end{align}
We now bound $B_2$.
\begin{align}
B_2 &= \left( \vxm^2 f(w)  \,\vnu , (1 - \rho)\nabla F(\vx_0) \, \right) - \left( \vxhm^2 f(\wh)  \,\vnuh , (1 - \rho)\nabla F(\vxh_0) \, \right)^h \nono \\
&= \left( \vxm^2 f(w)  \,\vnu, (1 - \rho)\left[\nabla F(\vx_0) - \nabla F(\vxh_0) \right] \right) + \left( \vxhm^2 \left[f(w) - f(\wh) \right] \vnuh , (1 - \rho)\nabla F(\vxh_0) \, \right) \nono \\
&\quad + \left( \left[\vxm^2 - \vxhm^2 \right] f(w)  \,\vnu + \vxhm^2 f(w)  \left[\vnu - \vnuh \right] , (1 - \rho)\nabla F(\vxh_0) \, \right) \nono \\
&\quad + \biggl[ \left( \vxhm^2 \left[ (I - I^h) f(\wh) \right] \vnuh, (1-\rho) \nabla F(\vxh_0) \right) + \left( \vxhm^2 \left[ (I^h - I) f(\wh) \right] \vnuh, (1-\rho) \nabla F(\vxh_0) \right)^h \biggr] \nono \\
&\quad + \biggl[\left( \vxhm^2 \,I^h(f(\wh))  \,\vnuh , (1 - \rho)\nabla F(\vxh_0) \, \right) - \left( \vxhm^2 \,I^h(f(\wh))  \,\vnuh , (1 - \rho)\nabla F(\vxh_0) \, \right)^h \biggr] =: \sum_{i=1}^5 B_{2,i}. \label{B2}
\end{align}
Using \eqref{nodegen}, \eqref{gFxgFxh} and \eqref{freg}, we have
\begin{align}
B_{2,1} &= \left( \vxm^2 f(w)  \,\vnu, (1 - \rho)\left[\nabla F(\vx_0) - \nabla F(\vxh_0) \right] \right) \leq C \, |f(w)|_0 \, |1-\rho|_0  \, |\thh|_{\zi} \leq C \, |\thh|_{0,\infty}. \label{B21}
\end{align}
Noting \eqref{sd_normF}, and using similar arguments to those used in proving \eqref{T21} and \eqref{T23}, we have
\begin{align}
B_{2,2} &= \left( \vxhm^2 \left[f(w) - f(\wh) \right] \vnuh , (1 - \rho)\nabla F(\vxh_0) \, \right) \leq C \, |f(w) - f(\wh)|_0 \, |1-\rho|_0 \leq C \left[h + |\zh|_0 \right], \label{B23}\\
B_{2,3} &=  \left( \left[\vxm^2 - \vxhm^2 \right] f(w)  \,\vnu + \vxhm^2 f(w)  \left[\vnu - \vnuh \right] , (1 - \rho)\nabla F(\vxh_0) \, \right) \nono \\
&\leq C \left[|\vx_\rho-\vxh_\rho|_0 + |\vnu - \vnuh|_0\right] |1-\rho|_0 \leq C \left[h + |\thh|_1 \right], \label{B22}\\
B_{2,4} &= \left( \vxhm^2 \left[ (I - I^h) f(\wh) \right] \vnuh, (1-\rho) \nabla F(\vxh_0) \right) + \left( \vxhm^2 \left[ (I^h - I) f(\wh) \right] \vnuh, (1-\rho) \nabla F(\vxh_0) \right)^h \nono \\
&\leq C \, |(I - I^h) f(\wh)|_{0} \, |1 - \rho|_0 \leq C\left[h + |\zh|_0 \right], \label{B24}\\
B_{2,5} &= \left( \vxhm^2 \,I^h(f(\wh))  \,\vnuh , (1 - \rho)\nabla F(\vxh_0) \, \right) - \left( \vxhm^2 \,I^h(f(\wh))  \,\vnuh , (1 - \rho)\nabla F(\vxh_0) \, \right)^h \nono \\
&\leq C \, h  \sum_{j=1}^J \left|I^h_j(f(\wh))\right|_{1,\sigma_j} \left| \vxhm^2 \, (1 - \rho) \, \inner{\nabla F(\vxh_0)}{\vnuh} \right|_{0,\sigma_j}  \leq C \, h \, |f(\wh)|_1 \, |1 - \rho|_0 \leq C \left[h + |\zh|_0 \right]. \label{B25}
\end{align}
We now bound $B_3$ in a similar way.
\begin{align}
B_3 &= \left( \vxhm^2 \left[\alpha \, \vxh_t + (1-\alpha) \innerb{\vxh_t}{\vnuh} \, \vnuh \right], (1 - \rho)\nabla F(\vxh_0) \, \right)^h \nono \\
& \qquad - \left( \vxm^2 \bigl[\alpha \,\vx_t + (1-\alpha) \innerb{\vx_t}{\vnu} \, \vnu \bigr], (1 - \rho)\nabla F(\vx_0) \, \right) \nono \\
&= \left( \vxm^2 \left[\alpha \,\vx_t + (1-\alpha) \innerb{\vx_t}{\vnu} \, \vnu \right], (1 - \rho) \left[\nabla F(\vxh_0) - \nabla F(\vx_0) \right] \right) \nono \\
&\quad + \left( \vxhm^2 \left[ \alpha \left(\vxh_t - I^h \vx_t \right) + (1 - \alpha)\innerb{( \vxh_t - I^h \vx_t)}{\vnuh} \, \vnuh \right], (1 - \rho)\nabla F(\vxh_0) \right)^h \nono \\
&\quad + \biggl[ \left( \left[ \vxhm^2 - \vxm^2 \right] \left[\alpha \,\vx_t + (1-\alpha) \innerb{\vx_t}{\vnu} \, \vnu \right], (1 - \rho)\nabla F(\vxh_0) \, \right) \nono \\
&\qquad \quad + (1 - \alpha) \left( \vxhm^2  \left[ \innerb{\vx_t}{(\vnuh - \vnu)} \, \vnu + \innerb{\vx_t}{\vnuh} \left(\vnuh - \vnu \right) \right], (1 - \rho)\nabla F(\vxh_0) \right) \biggr] \nono \\
&\quad + \left( \vxhm^2 \left[\alpha \, (I^h - I) \vx_t + (1-\alpha) \innerb{ (I^h - I) \vx_t}{\vnuh} \, \vnuh \right], (1 - \rho)\nabla F(\vxh_0) \right) \nono \\
&\quad + \biggl[\left( \vxhm^2 \left[\alpha \, I^h \vx_t + (1-\alpha) \innerb{I^h \vx_t}{\vnuh} \, \vnuh \right], (1 - \rho)\nabla F(\vxh_0) \, \right)^h \nono \\
& \qquad \quad - \left( \vxhm^2 \left[\alpha \, I^h \vx_t + (1-\alpha) \innerb{I^h \vx_t}{\vnuh} \, \vnuh \right], (1 - \rho)\nabla F(\vxh_0) \, \right) \biggr] =: \sum_{i=1}^5 B_{3,i}. \label{B3}
\end{align}
Using \eqref{nodegen}, \eqref{gFxgFxh} and \eqref{xh2}, we have
\begin{align}
B_{3,1} &= \left( \vxm^2 \left[\alpha \,\vx_t + (1-\alpha) \innerb{\vx_t}{\vnu} \, \vnu \right], (1 - \rho) \left[\nabla F(\vxh_0) - \nabla F(\vx_0) \right] \right) \nono \\
&\leq C \, |\vx_t|_0\,|1-\rho|_0\, |\thh|_{\zi} \leq C \, |\thh|_{0,\infty}. \label{B31}
\end{align}
Using \eqref{disnodegen}, \eqref{sd_normF} and \eqref{lumpev}, we have
\begin{align}
B_{3,2} &= \left( \vxhm^2 \left[ \alpha\left[\vxh_t - I^h \vx_t \right] + (1 - \alpha) \innerb{(\vxh_t - I^h \vx_t)}{\vnuh} \, \vnuh \right], (1 - \rho)\nabla F(\vxh_0) \right)^h \nono \\
&\leq C \, \|\thh_t\|_h \, \|1- \rho\|_h \leq C \, |\thh_t|_0 . \label{B33}
\end{align}
Noting \eqref{sd_normF} and using similar arguments to those used in proving \eqref{T11} and \eqref{T13} we have
\begin{align}
B_{3,3} &= \left( \left[\vxhm^2 - \vxm^2 \right] \left[\alpha \,\vx_t + (1-\alpha) \innerb{\vx_t}{\vnu} \, \vnu \right], (1 - \rho)\nabla F(\vxh_0) \, \right) \nono \\
&\quad + \, (1 - \alpha) \left( \vxhm^2 \left[ \innerb{\vx_t}{(\vnuh - \vnu)} \, \vnu + \innerb{\vx_t}{\vnuh} \left(\vnuh - \vnu \right) \right], (1 - \rho)\nabla F(\vxh_0) \right) \nono \\
&\leq C \left[| \vx_\rho-\vxh_\rho|_0+|\vnu - \vnuh|_0\right]|1 - \rho|_0 \leq C \left[ h + |\thh|_1 \right],\label{B32}\\
B_{3,4} &= \left( \vxhm^2 \left[\alpha \, (I^h - I) \vx_t + (1-\alpha) \innerb{(I^h - I) \vx_t}{\vnuh} \, \vnuh \right], (1 - \rho)\nabla F(\vxh_0) \right) \nono \\
& \leq C \, |(I - I^h) \vx_t|_0 \, |1 - \rho|_0 \leq C \, h,  \label{B34}\\
B_{3,5} &= \left( \vxhm^2 \left[\alpha \, I^h \vx_t + (1-\alpha) \innerb{I^h \vx_t}{\vnuh} \, \vnuh \right], (1 - \rho)\nabla F(\vxh_0) \, \right)^h \nono \\
& \qquad - \left( \vxhm^2 \left[\alpha \, I^h \vx_t + (1-\alpha) \innerb{I^h \vx_t}{\vnuh} \, \vnuh \right], (1 - \rho)\nabla F(\vxh_0) \, \right) \nono \\
&\leq C \, h  \sum_{j=1}^J \left|I^h_j(\vx_t)\right|_{1,\sigma_j} \left| \vxhm^2 (1 - \rho) \left[ \nabla F(\vxh_0) + \innerb{\nabla F(\vxh_0)}{\vnuh} \, \vnuh \right] \right|_{0,\sigma_j}
\leq C \, h \, |\vx_t|_1 \, |1 - \rho|_0 \leq C \, h. \label{B35}
\end{align}
Combining \eqref{bRHS} with \eqref{B1}--\eqref{B35}, we have
\begin{align*}
|b(0,t) - b^h(0,t)| \leq C \left[ h + |\thh|_{0,\infty} + |\thh_t|_0 + |\zh|_0 + |\thh|_1 \right].
\end{align*}
We remark that the above bound does not depend on $\rho$ and so also holds for $\rho=1$ and hence we have
\begin{align}
|b(0,t) - b^h(0,t)| + |b(1,t) - b^h(1,t)| \leq C \left[ h + |\thh|_{0,\infty} + |\thh_t|_0 + |\zh|_0 + |\thh|_1 \right]. \label{bRHSb}
\end{align}
Combining \eqref{bRHSb} with Sobolev embeddings, \eqref{xh2} and \eqref{gFxgFxh}, noting \eqref{interp_bc}, we have
\begin{align}
I_{3,2} &= \left[(b(\rho,t) - b^h(\rho,t)) (\vx_t\cdot (\nabla F(\vx) - \nabla F(\vxh)) ) \right]_0^1 \leq C |\thh|_{\zi} \left[ h + |\thh|_{0,\infty} + |\thh_t|_0 + |\zh|_0 + |\thh|_1  \right]. \label{T32}
\end{align}
Hence, combining \eqref{T3} with \eqref{T31} and \eqref{T32}, we have
\begin{align}
I_3 &\leq \frac{1}{2} \frac{d}{dt} \left[ \innerb{\vx_\rho}{\grad F(\vx)} \thh^T D^2 F(\vx) \, \thh \right]^1_0 + \frac{m^2 \alpha}{24} \, |\thh_t|_0^2 + C\left[ h^2 + |\thh|_{\zi}^2 + |\zh|_0^2  + |\thh|_1^2 + h^{-1} \, |\thh|_{\zi}^4 \right]. \label{T3b}
\end{align}
Combining \eqref{xLHSb}, \eqref{T1b}, \eqref{T2b} and \eqref{T3b}, we have
\begin{align}
\frac{1}{2} \frac{d}{dt} |\thh|_1^2 + \frac{m^2 \alpha}{8} \, |\thh_t|_0^2 &\leq \frac{1}{2} \frac{d}{dt} \left[ \innerb{\vx_\rho}{\grad F(\vx)} \thh^T D^2 F(\vx) \, \thh \right]^1_0 \nono \\
& \qquad + C \left[h^2 + |\thh|_{0,\infty}^2 + |\zh|_0^2 + |\thh|_1^2 + h^{-1} \, |\thh|_{0,\infty}^4 \right]. \label{xRHSb}
\end{align}
Multiplying \eqref{xRHSb} by $\egs$, for $\gamma \geq 1$, and integrating with respect to $s \in (0,t)$ with $t \leq \Ts$, and noting $|\thh(\cdot,0)| = 0$, we have
\begin{align}
& \frac{1}{2} \egt |\thh|_1^2 + \frac{\gamma}{2} \ints \egs |\thh|_1^2 \, \ds + \frac{m^2 \alpha}{8} \ints \egs |\thh_t|_0^2 \, \ds \nono \\
& \quad \leq \frac{1}{2} \egt \left[ \innerb{\vx_\rho}{\grad F(\vx)} \thh^T D^2 F(\vx) \, \thh \right]^1_0 + \frac{\gamma}{2} \ints \egs \left[\innerb{\vx_\rho}{\grad F(\vx)} \thh^T D^2 F(\vx) \, \thh \right]^1_0 \ds \nono \\
& \qquad + C \ints \egs \left[ h^2 + |\thh|_{\zi}^2 + |\zh|_0^2 + |\thh|_1^2 + h^{-1} \, |\thh|_{\zi}^4 \right] \ds \nono \\
& \quad =: I_4 + C \ints \egs \left[ h^2 + |\thh|_{\zi}^2 + |\zh|_0^2 + |\thh|_1^2 + h^{-1} \, |\thh|_{\zi}^4 \right] \ds. \label{Irhsx}
\end{align}
Using \eqref{Freg}, Sobolev embeddings, \eqref{xh2} and the inequality
\begin{align}
|\eta|_{\zi}^2  \leq C |\eta|_0 \|\eta\|_1  \leq \eps \, |\eta|_1^2 + C(\eps) |\eta|_0^2, \quad \mbox{for } \eta \in H^1(\Ip), \label{GNiq}
\end{align}
we see that
\begin{align}
I_4 &= \frac{1}{2} \egt \left[ \innerb{\vx_\rho}{\grad F(\vx)} \thh^T D^2 F(\vx) \, \thh \right]^1_0 + \frac{\gamma}{2}\ints \egs \left[ \innerb{\vx_\rho}{\grad F(\vx)} \thh^T D^2 F(\vx) \, \thh \right]^1_0 \ds \nono \\
& \leq \egt |\vx|_{\oi} |D^2 F(\vx)|_{\zi} |\thh|_{\zi}^2 + \gamma \ints \egs |\vx|_{\oi} |D^2 F(\vx)|_{\zi} |\thh|_{\zi}^2 \, \ds \nono\\
& \leq \frac{1}{4} \egt |\thh|_1^2 + C\egt  |\thh|_0^2+ \ints \egs \left[ \frac{\gamma}{4} \, |\thh|_1^2 + C\gamma  |\thh|_0^2 \right] \ds.  \label{I1}
\end{align}
Substituting \eqref{I1} into \eqref{Irhsx} and using \eqref{GNiq}, gives
\begin{align}
& \frac{1}{4} \egt |\thh|_1^2 + \frac{\gamma}{4} \ints \egs |\thh|_1^2 \, \ds + \frac{m^2 \alpha}{8} \ints \egs |\thh_t|_0^2 \, \ds \nono \\
&  \qquad\leq C \egt |\thh|_0^2 + \frac{C \gamma}2 \ints \egs |\thh|_{0}^2 \ds + C \ints \egs \left[ h^2 + |\zh|_0^2 + |\thh|_1^2 + h^{-1} \, |\thh|_{\zi}^4 \right] \ds. \label{Irhsx2}
\end{align}
Since $|\thh(\cdot,0)| = 0$, we have
\begin{align*}
\egt |\thh(\cdot,t)|_0^2 &= \ints \frac{d}{d s} \left( \egs |\thh|_0^2 \right) \ds \\
&\leq -\gamma \ints \egs |\thh|_0^2 \, \ds + 2 \ints \egs |\thh|_0 \, |\thh_t|_0 \, \ds \leq -\frac{\gamma}{2} \ints \egs |\thh|_0^2 \, \ds + \frac{2}{\gamma} \ints \egs |\thh_t|_0^2 \, \ds,
\end{align*}
and hence there exists $C_2\in \mathbb{R}_{>0}$ such that first two terms on the right hand side of \eqref{Irhsx2} can be bounded as follows
\begin{align}
C \egt |\thh|_0^2 + \frac{C\gamma}{2} \ints \egs |\thh|_0^2 \, \ds \leq \frac{2C_2}{\gamma} \ints \egs |\thh_t|_0^2 \, \ds. \label{infb}
\end{align}
Combining \eqref{Irhsx2} and \eqref{infb}, with $\gamma$ chosen large enough such that $\gamma\geq \max\{1,\frac{32C_2}{m^2 \alpha}\}$, yields the desired result. \hfill $\Box$

{\bf Proof of Lemma \ref{rde_lem}: }
In the proof of this lemma we follow the techniques used in \cite{BDS17}. 
We first note that from \eqref{cont} and \eqref{exp}, for $h \in (0,\hs]$ and $t \in [0,\Ts)$, we have
\begin{align}
\label{EZh_bound}
|\thh(\cdot,t)|_1^2 + |\zh(\cdot,t)|_0^2 \leq 2 C_1 h^2 e^{\gamma T} \leq 2 C_1 (\hs)^{\frac{1}{2}} h^{\frac{3}{2}} e^{\gamma \, T} \leq h^{\frac{3}{2}}.
\end{align}
Setting $\eta = \zh$ in \eqref{wwf}, subtracting the resulting equation from \eqref{wfem} with $\eta^h=\zh$, and noting \eqref{xrihxr_int}, gives
\begin{align*}
\left( \left(\vxm w \right)_t, \zh \right) &- \left( \left(\vxhm \, \wh \right)_t ,\zh \right)^h +  \left( \frac{1}{\vxhm}\zh_\rho, \zh_\rho \right)
+ \left( \psi \, w, \zh_\rho \right) - \left(\psi^h\,\wh,\zh_\rho^h \right)^h \\
&=   \left(w_\rho \left[\frac{1}{\vxhm} - \frac{1}{\vxm} \right]  , \zh_\rho \right)  + \left( \vxm \, g(v,w), \zh \right) - \left( \vxhm \,g(v^h,\wh), \zh \right)^h. 
\end{align*}
Since
\begin{align*}
\left( \left(\vxhm \, I^h w \right)_t ,\zh \right)^h -\left( \left(\vxhm \, \wh \right)_t ,\zh \right)^h & = \frac{1}{2} \frac{d}{dt} \left( \vxhm \, \zh ,\zh \right)^h
+\frac{1}{2} \left( \vxhm_t \, \zh, \zh \right)^h,
\end{align*}
we have
\begin{align}
&\frac{1}{2} \frac{d}{dt} \left( \vxhm \, \zh ,\zh \right)^h + \left(\frac{1}{\vxhm} \zh_\rho, \zh_\rho \right) \nono \\
& \quad = -\frac{1}{2} \left( \vxhm_t \, \zh, \zh \right)^h + \left[\left( \left(\vxhm \, I^h w \right)_t ,\zh \right)^h - \left( \left(\vxm \, w \right)_t ,\zh \right) \right] + \left(w_\rho \left[\frac{1}{\vxhm} - \frac{1}{\vxm} \right]  , \zh_\rho \right) \nono \\
&\qquad + \left[\left(\psi^h\,\wh,\zh_\rho \right)^h - \left( \psi \, w, \zh_\rho \right) \right] + \left[ \left( \vxm \, g(v,w), \zh \right) - \left( \vxhm\,g(v^h,\wh), \zh \right)^h \right] =: \sum_{i=1}^{5} T_i. \label{wRHS}
\end{align}
Using \eqref{disnodegen}, the left hand side of \eqref{wRHS} is bounded below by
\begin{align}
\frac{1}{2} \frac{d}{dt} \left( \vxhm \, \zh ,\zh \right)^h + \left(\frac{1}{\vxhm} \zh_\rho, \zh_\rho \right) \geq \frac{1}{2} \frac{d}{dt} \left( \vxhm \, \zh ,\zh \right)^h + \frac{1}{2 M} |\zh|_1^2. \label{wLHSb}
\end{align}
Now we bound $T_i$, $i=1,\dots,5$, by starting with $T_1$. Noting that $\vxhm_t = \inner{\vxh_{\rho,t}}{\vtauh}$ we have
\begin{align}
\label{T7}
T_1 = -\frac{1}{2} \left( \vxhm_t \, \zh, \zh \right)^h &= \frac{1}{2} \left( \innerb{\vxh_{\rho,t}}{(\vtau - \vtauh)} \, \zh, \zh \right) - \frac{1}{2} \left( \innerb{\vxh_{\rho,t}}{\vtau} \, \zh,\zh \right) \nono \\
& \qquad + \frac{1}{2} \left[\left(\innerb{\vxh_{\rho,t}}{\vtauh} \, \zh, \zh \right) - \left(\innerb{\vxh_{\rho,t}}{\vtauh} \, \zh, \zh \right)^h \right] =: \sum_{i=1}^3 T_{1,i}.
\end{align}
Using \eqref{xht0}, \eqref{tn0}, \eqref{EZh_bound} and \eqref{GNiq} gives
\begin{align}
T_{1,1} = \frac{1}{2} \left( \innerb{\vxh_{\rho,t}}{(\vtau - \vtauh)} \, \zh, \zh \right) &\leq \frac{1}{2} |\zh|_{\zi}^2 |\vxh_t|_1 |\vtau - \vtauh|_0 \leq C  |\zh|_{\zi}^2 \left[1 + h^{-1} |\thh_t|_0 \right]\left[h + |\thh|_1 \right]\nono\\
& \leq C h^{\frac{3}{4}} |\zh|_{\zi}^2 + C h^{\frac{1}{2}} \|\zh\|_1 |\thh_t|_0. \label{T71}
\end{align}
Using integration by parts with \eqref{xh2} and \eqref{xht0}, noting $\zh \in V^h_0$, yields
\begin{align}
\label{T72}
T_{1,2} &= - \frac{1}{2} \left( \innerb{\vxh_{\rho,t}}{\vtau} \, \zh,\zh \right) = \frac{1}{2} \left( \innerb{\vxh_t}{\vtau_\rho} \, \zh, \zh \right) +\left( \innerb{\vxh_t}{\vtau} \, \zh, \zh_\rho \right)\nono\\
&\leq C|\vxh_t|_0( |\vtau|_1  |\zh|_{\zi}^2 + |\zh|_{\zi} |\zh|_1) \leq C \left[1 + |\thh_t|_0\right] (|\zh|_{\zi}^2 +|\zh|_{\zi} |\zh|_1),
\end{align}
while \eqref{lumperr} and \eqref{xht0} yield
\begin{align}
\label{T74}
2 \, T_{1,3} &= \left(\innerb{\vxh_{\rho,t}}{\vtauh} \, \zh, \zh \right) - \left(\innerb{\vxh_{\rho,t}}{\vtauh} \, \zh, \zh \right)^h \nono \\
&\leq C \, h \, \sum_{j=1}^J |\zh|_{1,\sigma_j} |\innerb{\vxh_{\rho,t}}{\vtauh} \, \zh|_{0,\sigma_j} \leq C \, h \, |\zh|_{\zi} |\vxh_t|_1 |\zh|_1 \leq C |\zh|_{\zi} \left[h + |\thh_t|_0 \right] |\zh|_1.
\end{align}
Thus, noting that \eqref{xh2} and \eqref{whinf} imply that $|\zh|_{\zi} \leq C$, combining \eqref{T7} with \eqref{T71}--\eqref{T74}, and using \eqref{GNiq}, we have
\begin{align}
|T_1| \leq \frac{1}{8M} |\zh|_1^2 + \eps h \|\zh\|_1^2  + C \left[h^2 + |\zh|_{\zi}^2 + |\thh_t|_0^2 \right]
\leq \frac{1}{16M} |\zh|_1^2 + C \left[h^2 + |\thh_t|_0^2 + |\zh|_0^2 \right].
\label{T7b}
\end{align}
Noting $\vxm = \inner{\vx_{\rho}}{\vtau}$ and $\vxhm = \inner{\vxh_{\rho}}{\vtauh}$, in addition to $\inner{\vxh_\rho}{(\vtauh - \vtau)} = \frac{1}{2} \, \vxhm \, |\vtau - \vtauh|^2$ and $\inner{\vtauh}{\vtauh_t} = 0$, we have
\begin{align*}
\left(w  \left( \vxhm - \vxm \right)_t, \zh \right)&=\left( w \, ((\vxh_\rho - \vx_\rho)\cdot \vtau)_t, \zh \right) + \left( w \,  (\vxh_\rho\cdot (\vtauh - \vtau))_t, \zh \right) \nono\\
&=\left( w \, ((\vxh_\rho - \vx_\rho)\cdot \vtau)_t, \zh \right)+ \frac{1}{2} \left( w \left( \vxhm \, |\vtau - \vtauh|^2 \right)_t, \zh \right) \nono \\
&=\left( w \, ((\vxh_\rho - \vx_\rho)\cdot \vtau)_t, \zh \right)+ \frac{1}{2} \left( w \vxhm_t \, |\vtau - \vtauh|^2, \zh \right) + \frac{1}{2} \left( w \, \vxhm  (|\vtau - \vtauh|^2)_t, \zh \right) \nono\\
&=\left( w \, ((\vxh_\rho - \vx_\rho)\cdot \vtau)_t, \zh \right)+ \frac{1}{2} \left( w \, \vxhm_t \, |\vtau - \vtauh|^2, \zh \right) +  \left(w \, \vxhm \, \inner{(\vtau - \vtauh)}{(\vtau_t - \vtauh_t)}, \zh \right) \\
&=\left( w \, ((\vxh_\rho - \vx_\rho)\cdot \vtau)_t, \zh \right)+ \frac{1}{2} \left( w \, \vxhm_t \, |\vtau - \vtauh|^2, \zh \right) + \left( w \, \vxhm \, \inner{(\vtau - \vtauh)}{\vtau_t}, \zh \right) \\
&\quad - \left( w \, \vxhm \, \inner{\vtau}{\vtauh_t}, \zh \right)
\end{align*}
and hence
\begin{align}
T_2 &= \left(  \left(\vxhm \, I^h w \right)_t ,\zh \right)^h - \left(  \left(\vxm \, w \right)_t ,\zh \right) \nono \\
&= \left( w_t \left[ \vxhm - \vxm \right], \zh \right) + \left(w \left( \vxhm - \vxm \right)_t, \zh \right) + \left[\left( \vxhm \, (I^h - I)w_t, \zh \right) + \left( \vxhm_t \, (I^h - I)w, \zh \right) \right]\nono \\
& \quad + \left[ \left(  \left(\vxhm \, I^h w \right)_t ,\zh \right)^h - \left(  \left(\vxhm \, I^h w \right)_t ,\zh \right) \right] \nono\\
&= \left( w_t \left[ \vxhm - \vxm \right], \zh \right) + \left( w \, ((\vxh_\rho - \vx_\rho)\cdot \vtau)_t, \zh \right) + \frac{1}{2} \left( w \, \vxhm_t \, |\vtau - \vtauh|^2, \zh \right) +  \left( w \, \vxhm \, \inner{(\vtau - \vtauh)}{\vtau_t}, \zh \right)  \nono\\
&\quad - \left( w \, \vxhm \, \inner{\vtau}{\vtauh_t}, \zh \right)   + \left[\left( \vxhm \, (I^h - I)w_t, \zh \right) + \left( \vxhm_t \, (I^h - I)w, \zh \right) \right]\nono\\
&\quad + \left[ \left(  \left(\vxhm \, I^h w \right)_t ,\zh \right)^h - \left(  \left(\vxhm \, I^h w \right)_t ,\zh \right) \right]   =: \sum_{i=1}^7 T_{2,i}. \label{T8}
\end{align}
Using Sobolev embeddings, \eqref{xh2} and \eqref{xxh1}, we see that
\begin{align}
\label{T81}
T_{2,1} &= \left( w_t \left[ \vxhm - \vxm \right], \zh \right) \leq |w_t|_{\zi} |\vx - \vxh|_1 |\zh|_0 \leq C \left[ h + |\thh|_1 \right] |\zh|_0.
\end{align}
Applying integration by parts, \eqref{xh2}, Sobolev embeddings, \eqref{xtxht}, \eqref{xxh1} and noting that $\zh \in V_0^h$, we have that
\begin{align}
\label{T821}
T_{2,2} &= \left( w \, ((\vxh_\rho - \vx_\rho)\cdot\vtau)_t, \zh \right)
= \left( \inner{(\vx - \vxh)_t}{(\vtau \, w)_\rho}, \zh \right) + \left( w \, \innerb{(\vx - \vxh)_t}{\vtau}, \zh_\rho \right) + \left(w \, \inner{(\vxh_\rho - \vx_\rho)}{\vtau_t}, \zh \right) \nono \\
&\leq |(\vx - \vxh)_t|_0 \left[ \, |\vtau w|_1 \, |\zh|_{\zi} + |w|_{\zi} \, |\zh|_1 \right] + |w|_{\zi} \, |\vx - \vxh|_1 \, |\vtau_t|_0 \, |\zh|_{\zi} \nono \\
&\leq C \left[h + |\thh_t|_0 + |\thh|_1 \right] \left[ |\zh|_{\zi} + |\zh|_1 \right].
\end{align}
From \eqref{wreg}, \eqref{tn0},  \eqref{invh},  \eqref{hdef}, \eqref{xht0}, and \eqref{EZh_bound}, we have that
\begin{align}
\label{T8221}
T_{2,3} &= \frac{1}{2} \left( w \, \vxhm_t \, |\vtau - \vtauh|^2, \zh \right) \nono \\
&\leq \frac{1}{2} \, |w|_{\zi} \, |\vtau - \vtauh|_{\zi} \, |\vxh_t|_1 \, |\vtau - \vtauh|_0 \, |\zh|_{\zi} \nono \\
&\leq C \left[h^{\frac{1}{2}} + h^{-\frac{1}{2}} |\thh|_1 \right] \left[1 + h^{-1}|\thh_t|_0 \right] \left[h + |\thh|_1 \right] |\zh|_{\zi} \leq C  \left[ h + |\thh_t|_0 \right] |\zh|_{\zi} .
\end{align}
Using \eqref{disnodegen}, Sobolev embeddings, \eqref{xh2} and \eqref{tn0}, gives
\begin{align}
\label{T82221}
T_{2,4} &= \left( w \, \vxhm \, \inner{(\vtau - \vtauh)}{\vtau_t}, \zh \right) \leq C \, |w|_{\zi} \, |\vtau_t|_0 \, |\vtau - \vtauh|_0 \, |\zh|_{\zi} \leq C \left[ h + |\thh|_1 \right] |\zh|_{\zi}.
\end{align}
Setting $P^h := I - \vtauh \otimes \vtauh$, where $\otimes$ represents the outer product, and noting that $\vxhm \vtauh_t = P^h \vxh_{\rho,t}$, we have
\begin{align}
\label{T82222}
T_{2,5} &= - \left( w \, \vxhm \, \inner{\vtau}{\vtauh_t}, \zh \right) \nono \\
&= - \left( w \, \inner{P^h \vxh_{\rho,t}}{\vtau}, \zh \right) = \left(w \, \inner{P^h \thh_{\rho,t}}{\vtau}, \zh \right) - \left(w \, \inner{P^h I^h \vx_{\rho,t}}{\vtau}, \zh \right) =: T_{2,5,1} + T_{2,5,2}.
\end{align}
Since $P^h$ is constant on each sub-interval $\sigma_j$ and $\zh \in V_0^h$, using integration by parts over the sub-intervals yields
\begin{align}
T_{2,5,1} &= \left(w \, \inner{P^h \thh_{\rho,t}}{\vtau}, \zh \right)
= \sum_{j=1}^J \int_{\sigma_j} w \, \inner{P^h \thh_{\rho,t}}{\vtau} \, \zh \, \drho \nono \\
&= \sum_{j=1}^J \left[ w \, \inner{P^h \thh_t}{\vtau} \, \zh \right]_{\rho_{j-1}}^{\rho_j} - \sum_{j=1}^J \int_{\sigma_j} \inner{P^h \thh_t}{(w \, \vtau \, \zh)_\rho} \, \drho \nono \\
&= -\sum_{j=1}^{J-1} \left[ \inner{(P^h_{|_{\sigma_{j+1}}} - P^h_{|_{\sigma_j}})\thh_t(\rho_j,t)}{(w(\rho_j,t) \, \vtau(\rho_j,t) \, \zh(\rho_j,t))} \right] \nono \\
& \qquad \qquad - \left(\inner{P^h \thh_t}{(w \, \vtau)_\rho}, \zh \right) - \left(w \, \inner{P^h \thh_t}{\vtau}, \zh_\rho \right). \label{T822221}
\end{align}
To bound the first term in \eqref{T822221} we first note that
\begin{align*}
P^h_{|_{\sigma_{j+1}}} - P^h_{|_{\sigma_j}} &= \vtauh_{|_{\sigma_{j+1}}} \otimes ( \vtauh_{|_{\sigma_{j}}}  - \, \vtauh_{|_{\sigma_{j+1}}} )  + ( \vtauh_{|_{\sigma_{j}}}  - \, \vtauh_{|_{\sigma_{j+1}}} ) \otimes \vtauh_{|_{\sigma_{j}}} ,
\end{align*}
and
\begin{align*}
\vtauh_{|_{\sigma_{j+1}}} - \, \vtauh_{|_{\sigma_{j}}} &= \frac{1}{|(\vxh_\rho)_{|_{\sigma_{j+1}}}|} \left((\vxh_\rho)_{|_{\sigma_{j+1}}} - (\vxh_\rho)_{|_{\sigma_{j}}} \right) + \frac{\vtauh_{|_{\sigma_{j}}}}{|(\vxh_\rho)_{|_{\sigma_{j+1}}}|} \left( |(\vxh_\rho)_{|_{\sigma_{j}}}| - |(\vxh_\rho)_{|_{\sigma_{j+1}}}| \right).
\end{align*}
For any $\vec{\chi} \in \mathbb{R}^2$, we set $\vxi^h = \phi_j \vec{\chi}$, $j=1,\cdots,J-1$, in \eqref{xfem} to obtain
\begin{align*}
\inner{\left[ (\vxh_\rho)_{|_{\sigma_{j+1}}} - (\vxh_\rho)_{|_{\sigma_{j}}} \right]}{\vec{\chi}} &= - \left(h_{j+1} \, |(\vxh_\rho)_{|_{\sigma_{j+1}}}|^2 \, \innerb{\vnuh_{|_{\sigma_{j+1}}}}{\vec{\chi}} + h_{j} \, |(\vxh_\rho)_{|_{\sigma_{j}}}|^2 \, \innerb{\vnuh_{|_{\sigma_{j}}}}{\vec{\chi}}  \right) f(w^h(\rho_j,t)) \\
& \quad + \alpha \left( h_{j+1} \, |(\vxh_\rho)_{|_{\sigma_{j+1}}}|^2 + h_j \, |(\vxh_\rho)_{|_{\sigma_{j}}}|^2 \right) \innerb{\vxh_t(\rho_j,t)}{\vec{\chi}} \\
& \quad + (1 - \alpha) \, h_{j+1} \, |(\vxh_\rho)_{|_{\sigma_{j+1}}}|^2 \, \innerb{\vxh_t(\rho_j,t)}{\vnuh_{|_{\sigma_{j+1}}}} \, \innerb{\vnuh_{|_{\sigma_{j+1}}}}{\vec{\chi}} \\
& \qquad + (1 - \alpha) \, h_j \, |(\vxh_\rho)_{|_{\sigma_{j}}}|^2 \, \innerb{\vxh_t(\rho_j,t)}{\vnuh_{|_{\sigma_{j}}}} \, \innerb{\vnuh_{|_{\sigma_{j}}}}{\vec{\chi}}
\end{align*}
Combining the three equations above and using \eqref{disnodegen}, \eqref{freg}, \eqref{whinf} and \eqref{xh2}, we have
\begin{align}
|P^h_{|_{\sigma_{j+1}}} - P^h_{|_{\sigma_j}}| &\leq C \, | (\vxh_\rho)_{|_{\sigma_{j+1}}} - \, (\vxh_\rho)_{|_{\sigma_{j}}} | \leq C \, h \left[ |f(w^h)|_{\zi} + |\vxh_t(\rho_j,t)| \right] \leq C \, h \left[ 1 + |\thh_t(\rho_j,t)| \right].\label{dfg}
\end{align}
Hence, using \eqref{dfg}, \eqref{hdef}, \eqref{wreg}, \eqref{invh}, and \eqref{EZh_bound}, we have
\begin{align}
\label{T8222211}
\sum_{j=1}^{J-1} &\inner{(P^h_{|_{\sigma_{j+1}}} - P^h_{|_{\sigma_j}})\thh_t(\rho_j,t)}{(w(\rho_j,t) \, \vtau(\rho_j,t) \, \zh(\rho_j,t))} \nono \\
&\leq C \, h \sum_{j=1}^{J-1} \left[1 + |\thh_t(\rho_j,t)| \right] |\thh_t(\rho_j,t)| \, |w(\rho_j,t)| \, |\zh(\rho_j,t)| \nono \\
&\leq C \left[ |\thh_t|_0 + |\thh_t|_0^2 \right] |\zh|_{\zi}\leq C \, |\thh_t|_0 \, |\zh|_{\zi} + C \, h^{-\frac12} \,|\thh_t|_0^2 \, |\zh|_{0} \leq  C \, |\thh_t|_0 \, |\zh|_{\zi} + C \, |\thh_t|_0^2.
\end{align}
From \eqref{T8222211}, \eqref{disnodegen}, \eqref{xh2} and Sobolev embeddings, we have
\begin{align}
\label{T8222212}
T_{2,5,1} &= - \sum_{j=1}^{J-1} \inner{(P^h |_{\sigma_{j+1}} - P^h |_{\sigma_j})\thh_t(\rho_j,t)}{(w(\rho_j,t) \, \vtau(\rho_j,t) \, \zh(\rho_j,t))} -\left(\inner{P^h \thh_t}{(w \, \vtau)_\rho}, \zh \right) - \left(w \, \inner{P^h \thh_t}{\vtau}, \zh_\rho \right) \nono\\
& \leq C \, |\thh_t|_0 \, |\zh|_{\zi} + C \, |\thh_t|_0^2 + |P^h|_{\zi} \, |w \, \vtau |_1 \, |\thh_t|_0 \, |\zh|_{\zi} + |P^h|_{\zi} \, |w|_{\zi} \, |\thh_t|_0 \, |\zh|_1 \nono\\
& \leq C \, |\thh_t|_0 \, |\zh|_{\zi} + C \, |\thh_t|_0^2 + C \, |\thh_t|_0 \, |\zh|_1.
\end{align}
Since $P^h$ is symmetric and $P^h \vtau = \vtau - \vtauh + \frac{1}{2} |\vtau - \vtauh|^2 \vtauh$, using Sobolev embeddings, \eqref{lumpev}, \eqref{xh2}, \eqref{tn0} and the fact that $|\vtau - \vtauh| \leq |\vtau| + |\vtauh| \leq 2$, we have
\begin{align}
\label{T822222}
T_{2,5,2} &= - \left(w \, \inner{P^h I^h \vx_{\rho,t}}{\vtau}, \zh \right)
= -\left(w \, \inner{P^h \vtau}{I^h \vx_{\rho,t}}, \zh \right) \nono \\
&= \left(w \, \inner{(\vtauh - \vtau)}{I^h \vx_{\rho,t}}, \zh \right) - \frac{1}{2} \left(w \, |\vtau - \vtauh|^2 \innerb{\vtauh}{I^h \vx_{\rho,t}}, \zh \right) \nono \\
&\leq 2 \, |w|_{\zi} \, |\vtau - \vtauh|_0 \, |I^h \vx_t|_1 \, |\zh|_{\zi} \leq C \left[h + |\thh|_1 \right] |\zh|_{\zi}.
\end{align}
Thus, combining \eqref{T8}--\eqref{T822222}, we have
\begin{align}
\label{T82b}
T_{2,5} &\leq C \left[h + |\thh_t|_0 + |\thh|_1 \right] \left[ |\zh|_{\zi} + |\zh|_0 + |\zh|_1 \right] + C \, |\thh_t|_0^2.
\end{align}
Using \eqref{disnodegen}, \eqref{Ih}, Sobolev embeddings, \eqref{xh2} and \eqref{xht0}, we see that
\begin{align}
\label{T84}
T_{2,6} &= \left( \vxhm (I^h - I)w_t, \zh \right) + \left( \vxhm_t (I^h - I)w, \zh \right) \nono\\
&\leq C \, |(I - I^h)w_t|_0 \, |\zh|_0 +  |(I - I^h)w|_{\zi} \, |(\vxh_{\rho,t}\cdot \vtauh)|_0 \, |\zh|_0 \nono \\
&\leq C \, h \, |w_t|_1 \, |\zh|_0 + C \, h \, |w|_2 \, |\vxh_t|_1 \, |\zh|_0 \leq C \left[h + |\thh_t|_0 \right] |\zh|_0.
\end{align}
From (\ref{lumpev},b), Sobolev embeddings, \eqref{xh2}, \eqref{disnodegen} and \eqref{xht0}, we obtain
\begin{align}
\label{T85}
T_{2,7} &= \left(  \left(\vxhm \, I^h w \right)_t ,\zh \right)^h - \left(  \left(\vxhm \, I^h w \right)_t ,\zh \right) \nono \\
&\leq C \, h \sum_{j=1}^J |\zh|_{1,\sigma_j} \left| \left(\vxhm \, I^h w \right)_t \right|_{0,\sigma_j}
\leq C \, h \left[|I^h w|_{\zi} \, |\vxh_{t}|_1 + | I^h w_t|_0 \right] |\zh|_1 \leq C \left[h + |\thh_t|_0 \right] |\zh|_1.
\end{align}
Combining \eqref{T8} with \eqref{T81}--\eqref{T82221} and \eqref{T82b}--\eqref{T85}, we have
\begin{align*}
|T_2| \leq \frac{1}{8M} |\zh|_1^2 + C \left[ h^2 + |\zh|_{\zi}^2 + |\thh_t|_0^2 + |\zh|_0^2 + |\thh|_1^2 \right].
\end{align*}
and hence, using \eqref{GNiq}, we have
\begin{align}
\label{T8b}
|T_2| \leq \frac{1}{16M} |\zh|_1^2 + C \left[ h^2 + |\thh_t|_0^2 + |\zh|_0^2 + |\thh|_1^2 \right].
\end{align}
From Sobolev embeddings, \eqref{xh2} and \eqref{xrxhbr}, we gain
\begin{align}
T_3 &= \left(w_\rho \left[\frac{1}{\vxhm} - \frac{1}{\vxm} \right], \zh_\rho \right) \nono \\
&\leq |w|_{1,\infty} \left|\frac{1}{\vxm} - \frac{1}{\vxhm} \right|_0 |\zh|_1
\leq C \left[h + |\thh|_1 \right] |\zh|_1 \leq \frac{1}{16 M} |\zh|_1^2 + C \left[h^2 + |\thh|_1^2 \right]. \label{T4}
\end{align}
We now bound $T_4$.
\begin{align}
T_4 &= \left(\psi^h\,\wh,\zh_\rho \right)^h - \left( \psi \, w, \zh_\rho \right) \nono \\
&= \left(\left[ \psi^h - I^h \psi \right] w^h + I^h\psi \left[ w^h - I^h w \right], \zh_\rho \right)^h + \biggl[ \left(I^h\psi \, I^hw, \zh_\rho \right)^h - \left(I^h\psi \, I^hw, \zh_\rho \right)  \nono \\
&\quad + \left( (I^h - I)\psi\,  I^hw, \zh_\rho \right) + \left( \psi \, (I^h - I) w , \zh_\rho \right) \biggr] =: T_{4,1} + T_{4,2}. \label{T5}
\end{align}
From \eqref{whinf}, Sobolev embeddings, \eqref{lumpev}, \eqref{xh2}, \eqref{Ih} and \eqref{psipsih0}, we have
\begin{align}
T_{4,1} &= \left(\left[ \psi^h - I^h \psi \right] w^h + I^h\psi \left[ w^h - I^h w \right], \zh_\rho \right)^h \nono \\
&\leq \max\{|w^h|_{0,\infty},|I^h \psi|_{\zi}\} \left[\|(I - I^h)\psi\|_h + \|\psi - \psi^h\|_h  + \|\zh\|_h \right] \|\zh_\rho\|_h \nono \\
&\leq C \left[h + |\thh_t|_0 + |\zh|_0 + |\thh|_1 \right] |\zh|_1, \label{T51}
\end{align}
while using Sobolev embeddings together with (\ref{lumpev},b), \eqref{xh2} and \eqref{Ih} gives
\begin{align}
T_{4,2} &= \left(I^h\psi \, I^hw, \zh_\rho \right)^h - \left(I^h\psi \, I^hw, \zh_\rho \right)  + \left( (I^h - I) \psi \, I^hw, \zh_\rho \right) + \left( \psi \,(I^h - I) w, \zh_\rho \right)  \nono \\
&\leq C \, h \, \sum_{j=1}^J |I^h_j(\psi) I^h_j(w)|_{1,\sigma_j} |\zh_\rho|_{0,\sigma_j} + \max\{|I^h w|_{0,\infty},|\psi|_{0,\infty}\} \left[|(I - I^h)\psi|_0 + |(I - I^h) w|_0 \right] |\zh|_1\nono\\
& \leq C \, h \left[ |I^h\psi \, I^hw|_1 + |\psi|_1 \, |w|_1 \right] |\zh|_1 \leq C \, h \, |\zh|_1. \label{T53}
\end{align}
Combining \eqref{T5}--\eqref{T53}, we have
\begin{align}
|T_4| &\leq \frac{1}{16 M} |\zh|_1^2 + C \left[h^2 + |\thh_t|_0^2 + |\zh|_0^2  + |\thh|_1^2 \right]. \label{T5b}
\end{align}
By using the continuity of $g$, we bound $T_5$ in the following way
\begin{align}
T_5 &= \left( \vxm \, g(v,w), \zh \right) - \left( \vxhm \,g(v^h,\wh), \zh \right)^h \nono \\
&= \left(\left[\vxm - \vxhm \right] g(v,w), \zh \right) + \left [ \left( \vxhm (I - I^h)g(v,w), \zh \right) + \left( \vxhm I^h g(v,w), \zh \right) - \left( \vxhm I^h g(v,w), \zh \right)^h \right] \nono \\
&\quad + \left(\vxhm \left( \left[ g(I^h v, I^h w) - g(I^h v, \wh) \right] + \left[ g(I^h v, \wh) - g(v^h, \wh) \right] \right), \zh \right)^h =: \sum_{i=1}^3 T_{5,i}. \label{T6}
\end{align}
Using \eqref{xreg}, \eqref{wreg}, \eqref{greg} and \eqref{xxh1}, gives
\begin{align}
T_{5,1} &= \left(\left[\vxm - \vxhm \right] g(v,w), \zh \right) \leq |g(v,w)|_{0,\infty} \, | \vx - \vxh |_1 \, |\zh|_0 \leq C \left[h + |\thh|_1 \right] |\zh|_0. \label{T61}
\end{align}
From \eqref{disnodegen}, (\ref{lumpev},b), \eqref{Ih}, \eqref{xreg}, \eqref{wreg} and \eqref{greg}, we have
\begin{align}
T_{5,2} &= \left( \vxhm (I - I^h)g(v,w), \zh \right) + \left( \vxhm I^h g(v,w), \zh \right) - \left( \vxhm I^h g(v,w), \zh \right)^h \nono \\
&\leq C \, h \, |g(v,w)|_1 |\zh|_0 + C \, h \sum_{j=1}^J \left|I^h_j(g(v,w))\right|_{1,\sigma_j} \, \left|\vxhm \, \zh \right|_{0,\sigma_j} \nono \\
&  \leq C \, h \, |g(v,w)|_1 |\zh|_0 \leq C \, h \, |g'(v,w)|_{\zi} \left[|v|_1 + |w|_1 \right] |\zh|_0 \leq C \, h \, |\zh|_0. \label{T63}
\end{align}
Using \eqref{disnodegen}, \eqref{greg}, \eqref{xreg}, \eqref{whinf}, \eqref{lumpev}, \eqref{Ih}, \eqref{xh2} and \eqref{psipsih0}, we have
\begin{align}
T_{5,3} &= \left(\vxhm \left( \left[ g(I^h v, I^h w) - g(I^h v, \wh) \right] + \left[ g(I^h v, \wh) - g(v^h, \wh) \right] \right), \zh \right)^h \nono \\
&\leq C \left[ \|\zh\|_h + \|I^h v - v^h\|_h \right] \|\zh\|_h \leq C \left[ h + |\thh_t|_0 + |\zh|_0 + |\thh|_1 \right] |\zh|_0. \label{T64}
\end{align}
Thus, combining \eqref{T6}--\eqref{T64}, yields
\begin{align}
|T_5| & \leq C \left[h^2 + |\thh_t|_0^2 + |\zh|_0^2 + |\thh|_1^2 \right]. \label{T6b}
\end{align}
We now combine \eqref{wRHS}, \eqref{wLHSb}, \eqref{T7b}, \eqref{T8b}, \eqref{T4}, \eqref{T5b} and \eqref{T6b} to obtain
\begin{align}
&\frac{1}{2} \frac{d}{dt} \left( \vxhm \, \zh ,\zh \right)^h + \frac{1}{4 M} |\zh|_1^2 \leq  C \left[h^2 + |\thh_t|_0^2 + |\zh|_0^2 + |\thh|_1^2 \right]. \label{wRHSb}
\end{align}
Multiplying \eqref{wRHSb} by $\egs$, for $\gamma \geq 1$, integrating with respect to $s \in (0,t)$, with $t \leq \Ts$,  and noting $|\zh(\cdot,0)| = 0$, we have
\begin{align}
\frac{1}{2} \egt \left( \vxhm \, \zh ,\zh \right)^h + \frac{\gamma}{2} \ints \egs  \left( \vxhm \, \zh ,\zh \right)^h \ds &+ \frac{1}{4 M} \ints \egs |\zh|_1^2 \, \ds \nono \\
&  \leq C_3 \ints \egs \left[h^2 + |\thh_t|_0^2 + |\zh|_0^2  + |\thh|_1^2 \right] \ds. \label{Irhsw}
\end{align}
From \eqref{disnodegen} and \eqref{lumpev}, we have
\begin{align*}
\frac{1}{2} \egt \left( \vxhm \zh, \zh\right)^h + \frac{\gamma}{2} \ints \egs \left( \vxhm \zh, \zh\right)^h \ds  \geq \frac{m}{4} \egt |\zh|_0^2 + \frac{\gamma \, m}{4} \ints \egs |\zh|_0^2 \, \ds,
\end{align*}
which together with (\ref{Irhsw}) yields the desired result. \hfill $\Box$

{\bf Proof of Theorem \ref{thmH1fbrde}:}
Multiplying \eqref{c3} by $\omega$, where $\omega \in \mathbb{R}_{>0}$ is chosen such that $C_3 \omega \leq \frac{m^2 \alpha}{32}$, and adding the resulting inequality to (\ref{sl}), for $t \in [0,\Ts)$, we have
\begin{align*}
\frac{1}{4} \egt |\thh|_1^2 + \frac{m \, \omega}{4} \egt |\zh|_0^2 &+ \frac{m^2 \alpha}{32}  \ints \egs |\thh_t|_0^2 \, \ds + \frac{\omega}{4 M} \ints \egs |\zh|_1^2 \, \ds \nono \\
& \qquad \leq C(1+\omega) \ints \egs \left[h^2 + |\thh|_1^2 + |\zh|_0^2 + h^{-1}|\thh|_{\zi}^4 \right] \ds.
\end{align*}
An application of Gronwall's lemma then gives
\begin{align*}
\sup_{s \in [0,\Ts]} \egs \left[\frac{1}{4}  |\thh|_1^2 + \frac{m \, \omega}{4} |\zh|_0^2\right] &+ \frac{m^2 \alpha}{32} \intsT \egs |\thh_t|_0^2 \, \ds + \frac{\omega}{4 M} \intsT \egs |\zh|_1^2 \, \ds \\
& \qquad \leq C_{\omega,\gamma} \intsT \egs \left[ h^2 + h^{-1} |\thh|_{\zi}^4 \right] \ds
\end{align*}
where $C_{\omega,\gamma}$ depends on $\omega$, $\gamma$ and $T$, but not $\Ts$. Dividing by $\tilde{C} = \min\{ \frac{1}{4}, \frac{m \, \omega}{4} , \frac{m^2 \alpha}{32}, \frac{\omega}{4 M} \}$, we obtain
\begin{align}
\sup_{s \in [0,\Ts]} \egs \left[ |\thh|_1^2 + |\zh|_0^2 \right] + \intsT \egs \left(|\thh_t|_0^2+|\zh|_1^2\right) \ds \leq C_1 h^2 + C h^{-1} \intsT \egs |\thh|_{\zi}^4 \, \ds. \label{pre_thm}
\end{align}
Using \eqref{GNiq}, \eqref{infb} and \eqref{cont}, for $t \in [0,\Ts)$, we have
\begin{align*}
\egt |\thh(\cdot,t)|_{0,\infty}^2 &\leq \egt |\thh(\cdot,t)|_1^2 + C \egt |\thh(\cdot,t)|_0^2 \leq C C_1 h^2,
\end{align*}
and hence, for $t \in [0,\Ts)$, we have
\begin{equation*}
C h^{-1} \ints \egs |\thh|_{\zi}^4 \, \ds \leq C h^{-1} e^{\gamma T} \ints \left( \egs |\thh|_{\zi}^2 \right)^2 \ds \leq C (C_1)^2 T e^{\gamma T} h^3
\end{equation*}
which, together with (\ref{pre_thm}), and on noting \eqref{exp}, yields
\begin{align}
\label{pre_cont}
& \sup_{s \in [0,\Ts]} \egs \left[ |\thh|_1^2 + |\zh|_0^2 \right] + \intsT \egs \left(|\thh_t|_0^2+|\zh|_1^2\right) \ds \nono \\
& \qquad\qquad \leq C_1 h^2 +C (C_1)^2 T e^{\gamma T} h^3 \leq C_1 h^2 + C C_1 T (\hs)^{\frac{1}{2}} h^{2} \leq C_1 h^2 + \frac{1}{2} C_1 h^2 \leq \frac{3}{2} C_1 h^2,
\end{align}
provided $\hs$ is chosen small enough. We now follow the argument in \cite{DD09} to show that $\Ts = T$. If it were not the case that $\Ts=T$ we would have $\Ts < T$, and using \eqref{nodegen}, \eqref{Ihinf}, \eqref{invh}, \eqref{hdef}, \eqref{xh2}, \eqref{pre_cont} and \eqref{exp}, for $\rho \in [0,1]$, we would have
\begin{align*}
|\vxh_\rho(\rho,\Ts)| &\leq |\vx_\rho(\rho,\Ts)| + |\vx_\rho(\rho,\Ts) - \vxh_\rho(\rho,\Ts)| \\
&\leq M + |(I - I^h)\vx(\rho,\Ts)|_{\oi} + |\thh|_{\oi} \\
&\leq M + C \, h^{\frac{1}{2}} |\vx|_2 + C \, h^{-\frac{1}{2}} |\thh|_1 \leq M + C \, h^{\frac{1}{2}} \left[ |\vx|_2 + e^{\frac{\gamma}{2} T} \right] \leq M+C\beta \leq \frac{3}{2}M,
\end{align*}
provided that $\beta$ is chosen small enough, and similarly
\begin{align*}
|\vxh_\rho(\rho,\Ts)| &\geq |\vx_\rho(\rho,\Ts)| - |\vx_\rho(\rho,\Ts) - \vxh_\rho(\rho,\Ts)| \geq \frac{3m}{4}.
\end{align*}
Using \eqref{Ihinf}, \eqref{hdef}, \eqref{wreg}, \eqref{pre_cont}, \eqref{exp} and Sobolev embeddings, we would also gain
\begin{align*}
|w^h(\cdot,\Ts)|_{\zi} &\leq |I^hw(\cdot,\Ts)|_{\zi} + |I^hw(\cdot,\Ts) - w^h(\cdot,\Ts)|_{\zi} \\
&\leq |I^hw(\cdot,\Ts)|_{\zi} + C h^{-\frac12}|\zh(\cdot,\Ts)|_{0} \\
&\leq \|w\|_{C([0,T];L^{\infty}(\Ip))} + C h^{\frac{1}{2}} e^{\frac{\gamma}{2}T}  \\
&\leq(C_w+C\beta) \|w\|_{C([0,T];H^1(\Ip))}  \leq \frac{3}{2} C_w \|w\|_{C([0,T];H^1(\Ip))},
\end{align*}
provided that $\beta$ is chosen small enough. Thus we could then extend the discrete solution to an interval $[0,\Ts + \delta]$ for some $\delta > 0$ with
\begin{align*}
\frac{m}{2} \leq |\vxh_\rho| \leq 2 M \qquad \text{ in } [0,1] \times [0,\Ts + \delta] \\
\|w^h\|_{C([0,\Ts + \delta]; L^\infty(\Ip))} \leq 2 C_w \|w\|_{C([0,T]; H^1(\Ip))} \\
\sup_{s \in [0,\Ts + \delta]} \egs \left[ |\thh|_1^2 + |\zh|_0^2 \right] + \int_{0}^{\Ts + \delta} \egs \left[ |\thh_t|_0^2 + |\zh|_1^2 \right] \, \ds < 2 C_1 h^2
\end{align*}
which contradicts the definition of $\Ts$. Therefore $\Ts = T$ and from \eqref{pre_cont}, \eqref{Ih}, \eqref{xreg} and \eqref{wreg} we obtain the desired result. \hfill $\Box$

\section{Numerical results}
We investigate the experimental order of convergence (eoc) of a fully discrete finite element approximation of \eqref{xfem} and \eqref{wfem} and then we conclude with two simulations of diffusion induced grain boundary motion.

Let $0 = t_0 < t_1 < \cdots < t_{N-1} < t_N = T$ be a uniform partition of $[0,T]$, with $\Delta t \, N = T$. We denote $\vX^n\in [V^h]^2$ to be the fully discrete finite element approximation of $\vx(\cdot,t^n)$, and denote $W^n\in V^h$ to be the fully discrete finite element approximation of $w(\cdot,t^n)$. We set $\vX^{0}(\cdot) := I^h \vx^{\, 0}(\cdot)$ and $W^0(\cdot) = I^h w^0(\cdot)$ in $[0,1]$, such that $\vx^{\, 0}$ and $w^0$ satisfy their respective boundary conditions \eqref{Feq}--\eqref{xFp} and \eqref{wbc}, and we define the discrete time derivative to be
\begin{align*}
D_t a^n := \frac{a^n - a^{n-1}}{\Delta t_n}. 
\end{align*}
On each $\sigma_j, \, j=1,\dots,J$, we assign to each element $\vX^n \in [V^h]^2$ a piecewise constant discrete unit tangent and normal, denoted by $\vTau^n$ and $\vNu^n$, approximating $\vtau(\cdot,t^n)$ and $\vnu(\cdot,t^n)$ respectively, and a piecewise linear tangential and normal velocity, denoted by $\Psi^n$ and $V^n$, approximating $\psi(\cdot,t^n)$ and $v(\cdot,t^n)$ respectively, which take the form
\begin{align*}
\vTau^n := \frac{\vX^n_\rho}{\vXmn{n}}, \quad \vNu^n := (\vTau^n)^\perp, \quad \Psi^n := D_t \vX^n \cdot \vTau^n, \quad V^n := D_t \vX^n \cdot \vNu^n, \qquad \text{ on } \sigma_j, \, j=1,\dots,J.
\end{align*}
The fully discrete finite element form of \eqref{xfem} takes the form: given $\vX^{n-1} \in [V^h]^2$ and $W^{n-1} - \gd \in V^h_0$, find $\vX^n \in [V^h]^2$ such that
\begin{align}
& \left( \vXmn{n-1}^2 \left[ \alpha \, D_t \vX^n + (1 - \alpha) \innerbb{D_t \vX^n}{\vNu^{n-1}} \vNu^{n-1} \right], \vxi^h \, \right)^h + \left( \vX^n_\rho, \vxi^h_\rho \right) \nono \\
& \qquad = \left[ \innerbb{\vX^{n}_\rho}{\nabla F(\vX^{n-1})} \innerbb{\vxi^h}{\nabla F(\vX^{n-1})} \right]^1_0 + \left(\vXmn{n-1}^2 f(W^{n-1}) \, \vNu^{n-1}, \vxi^h \right)^h \qquad \forall \, \vxi^h \in [V^h]^2 \label{xfea}
\end{align}
with the additional boundary constraint
\begin{align}
\label{dis_bc}
\inner{D_t \vX^n(\rho)}{\grad F(\vX^{n-1}(\rho))} = 0 \qquad \text{ for } \rho \in \{0,1\}.
\end{align}
The fully discrete finite element approximation of \eqref{wfem} takes the form: given $\vX^{n-1}, \vX^n \in [V^h]^2$ and $W^{n-1} - \gd \in V_0^h$, find $W^n - \gd \in V_0^h$ such that
\begin{align}
\label{wfea}
 \left( D_t \left[ \vXmn{n} \, W^n \right], \eta^h \right)^h  + \left( \frac{W^n_\rho}{\vXmn{n}}, \eta^h_\rho \right) + \left( \Psi^{n} \, W^{n}, \eta^h_\rho \right)^h = \left( \vXmn{n} \, g(V^n,W^{n-1}), \eta^h \right)^h \qquad \forall \, \eta^h \in V_0^h.
\end{align}
Setting $\alpha = 1$ in \eqref{xfea}--\eqref{dis_bc}  yields the discretisation is presented in \cite{DE98}, while neglecting the boundary condition, \eqref{dis_bc}, yields the discretisation presented in \cite{BDS17} for the closed curve configuration. For all our computations we use a uniform mesh with $h \, J = 1$.

\subsection{Estimated order of convergence}

We monitor the following error estimates
\begin{align*}
\Er{1} := \sup_{n = 0,\dots,N} |\vE^n|_1^2, \quad \Er{2} := \sum_{n=1}^N \Dt \, |D_t \vE^n|_0^2, \quad \Er{3} := \sup_{n = 0,\dots,N} |Z^n|_0^2, \quad \Er{4} := \sum_{n=1}^n \Dt \, |Z^n|_1^2
\end{align*}
where $\vE^n := I^h\vx^n - \vX^n$ and $Z^n := I^hw^n - W^n$, and quantify them using the estimated order of convergence (eoc)
\begin{align*}
\SwapAboveDisplaySkip
eoc_{i,j} := \frac{\ln(\Er{i,j+1}) - \ln(\Er{i,j})}{\ln(h_{j+1}) - \ln(h_j)},~~i=1,\ldots,4,~j=1,\cdots,5
\end{align*}
where $i$ corresponds to the error $\Er{i}$, $j$ corresponds to the relative mesh size $h_j$ and $\Er{i,j}$ corresponds to the error $\Er{i}$ at mesh level $j$.

{\bf Example 1}

Considering $\Omega := \mathbb{R} \times \mathbb{R}_+$ and taking $\Gamma(0)$ to be a semi-circle with radius one centred around the origin as well as $\gd = 0$, we see that the solution to \eqref{xeq}--\eqref{wbc}, for $(\rho,t) \in [0,1] \times [0,0.4]$, is
\begin{align*}
\vx(\rho,t) = \sqrt{1 - t} \left( \cos(\pi \rho), \sin(\pi \rho) \right)^T, \quad w(\rho,t) = (1 - t) \sin(\pi \rho)
\end{align*}
with the data
\begin{align*}
\SwapAboveDisplaySkip
f(w) = -\frac{w^2}{2(1-t)^{\frac{5}{2}}} - \frac{\cos^2(\pi \rho)}{2\sqrt{1-t}}, \quad g(v,w) = -\frac{w}{2(1-t)}.
\end{align*}
In Tables \ref{tab:a1dth2} and \ref{tab:az1dth2} we demonstrate the experimental convergence for $\Dt = h^2$ with $\alpha = 1$ and $\alpha = 0.1$ respectively and in Tables \ref{tab:a1dth} and \ref{tab:az1dth} we demonstrate the experimental convergence for $\Dt = 0.4 h$ with $\alpha = 1$ and $\alpha = 0.1$ respectively. These eoc results we observe in Tables \ref{tab:a1dth2}--\ref{tab:az1dth} are consistent with the eoc results presented in \cite{BDS17}, in particular we see first order convergence for $\Er{i}$, $i=1,\ldots,4$ with $\Dt = 0.4 h$, and second order convergence for $\Er{i}$, $i=1,\ldots,4$ with $\Dt = h^2$. We also observe that for $\Dt = h^2$, there is a sizeable reduction in the magnitudes of $\Er{1}$ and $\Er{2}$ when $\alpha$ is reduced from $\alpha=1$ to $\alpha=0.1$. We highlight here the different parameters $J$ and $N$ used for Tables \ref{tab:a1dth2} and \ref{tab:az1dth2} compared to Tables \ref{tab:a1dth} and \ref{tab:az1dth}.

\begin{table}[ht!]
\centering
\begin{tabular}{|c|c|c|c|c|c|c|c|c|c|}
\hline
$J$ & $N$ & $\Er{1} \times e^{-4}$ & $eoc_1$ & $\Er{2} \times e^{-5}$ & $eoc_2$ & $\Er{3} \times e^{-5}$ & $eoc_3$ & $\Er{4} \times e^{-5}$ & $eoc_4$ \\
\hline
\hline
10 & 80 & 44.54 & - & 147.0 & - & 1.123 & - & 5.522 & - \\
\hline
20 & 320 & 5.587 & \ter{3.55} & 13.34 & \ter{3.46} & 0.06858 & \ter{4.03} & 0.3491 & \ter{3.98} \\
\hline
40 & 1280 & 0.3812 & \ter{3.88} & 0.9244 & \ter{3.85} & 0.004296 & \ter{4.00} & 0.02186 & \ter{4.00} \\
\hline
80 & 5120 & 0.02436 & \ter{3.97} & 0.05933 & \ter{3.96} & 0.0002686 & \ter{4.00} & 0.001367 & \ter{4.00} \\
\hline
160 & 20480 & 0.00153 & \ter{3.99} & 0.003733 & \ter{3.99} & 0.00001679 & \ter{4.00} & 0.00008549 & \ter{4.00} \\
\hline
\end{tabular}
\caption{$\alpha = 1$, $\Dt = h^2$. }
\label{tab:a1dth2}
\end{table}

\begin{table}[ht!]
\centering
\begin{tabular}{|c|c|c|c|c|c|c|c|c|c|}
\hline
$J$ & $N$ & $\Er{1} \times e^{-4}$ & $eoc_1$ & $\Er{2} \times e^{-5}$ & $eoc_2$ & $\Er{3} \times e^{-5}$ & $eoc_3$ & $\Er{4} \times e^{-4}$ & $eoc_4$ \\
\hline
\hline
10 & 80 & 2.904 & - & 8.342 & - & 2.415 & - & 1.189 & - \\
\hline
20 & 320 & 0.1855 & \ter{3.97} & 0.6048 & \ter{3.79} & 0.1519 & \ter{4.00} & 0.07460 & \ter{3.99} \\
\hline
40 & 1280 & 0.01166 & \ter{3.99} & 0.03941 & \ter{3.94} & 0.009504 & \ter{4.00} & 0.004667 & \ter{4.00} \\
\hline
80 & 5120 & 0.0007296 & \ter{4.00} & 0.002490 & \ter{3.98} & 0.0005942 & \ter{4.00} & 0.0002918 & \ter{4.00} \\
\hline
160 & 20480 & 0.00004562 & \ter{4.00} & 0.0001560 & \ter{4.00} & 0.00003714 & \ter{4.00} & 0.00001824 & \ter{4.00} \\
\hline
\end{tabular}
\caption{$\alpha = 0.1$, $\Dt = h^2$. }
\label{tab:az1dth2}
\end{table}

\begin{table}[ht!]
\centering
\begin{tabular}{|c|c|c|c|c|c|c|c|c|c|}
\hline
$J$ & $N$ & $\Er{1}\times e^{-3}$ & $eoc_1$ & $\Er{2} \times e^{-3}$ & $eoc_2$ & $\Er{3} \times e^{-6}$ & $eoc_3$ & $\Er{4} \times e^{-7}$ & $eoc_4$ \\
\hline
\hline
40 & 80 & 7.651 & - & 2.111 & - & 1.608 & - & 9.976 & - \\
\hline
80 & 160 & 2.325 & \ter{1.72} & 0.6703 & \ter{1.66} & 0.3149 & \ter{2.35} & 1.976 & \ter{2.34} \\
\hline
160 & 320 & 0.6454 & \ter{1.85} & 0.1909 & \ter{1.81} & 0.06636 & \ter{2.25} & 0.4782 & \ter{2.05} \\
\hline
320 & 640 & 0.1704 & \ter{1.92} & 0.05110 & \ter{1.90} & 0.01494 & \ter{2.15} & 0.1211 & \ter{1.98} \\
\hline
640 & 1280 & 0.04379 & \ter{1.96} & 0.01323 & \ter{1.95} & 0.003523 & \ter{2.08} & 0.03073 & \ter{1.98} \\
\hline
\end{tabular}
\caption{$\alpha = 1$, $\Dt = 0.4 h$. }
\label{tab:a1dth}
\end{table}

\begin{table}[ht!]
\centering
\begin{tabular}{|c|c|c|c|c|c|c|c|c|c|}
\hline
$J$ & $N$ & $\Er{1} \times e^{-3}$ & $eoc_1$ & $\Er{2} \times e^{-3}$ & $eoc_2$ & $\Er{3}\times e^{-6}$ & $eoc_3$ & $\Er{4} \times e^{-7}$ & $eoc_4$ \\
\hline
\hline
40 & 80 & 6.874 & - & 1.931 & - & 1.591 & - & 10.41 & - \\
\hline
80 & 160 & 2.205 & \ter{1.64} & 0.6407 & \ter{1.59} & 0.3148 & \ter{2.34} & 1.856 & \ter{2.49} \\
\hline
160 & 320 & 0.6285 & \ter{1.81} & 0.1866 & \ter{1.78} & 0.06678 & \ter{2.24} & 0.4542 & \ter{2.03} \\
\hline
320 & 640 & 0.1681 & \ter{1.90} & 0.05053 & \ter{1.88} & 0.01509 & \ter{2.15} & 0.1182 & \ter{1.94} \\
\hline
640 & 1280 & 0.04351 & \ter{1.95} & 0.01316 & \ter{1.94} & 0.003563 & \ter{2.08} & 0.03052 & \ter{1.95} \\
\hline
\end{tabular}
\caption{$\alpha = 0.1$, $\Dt = 0.4 h$. }
\label{tab:az1dth}
\end{table}

\subsection{Numerical simulations of diffusion induced grain boundary motion}
We conclude the numerical results with two simulations of diffusion induced grain boundary motion (DIGM).

{\bf Example 2}

\begin{figure}[ht!]
\centering
\subfigure
{ \includegraphics[width=0.45\textwidth]{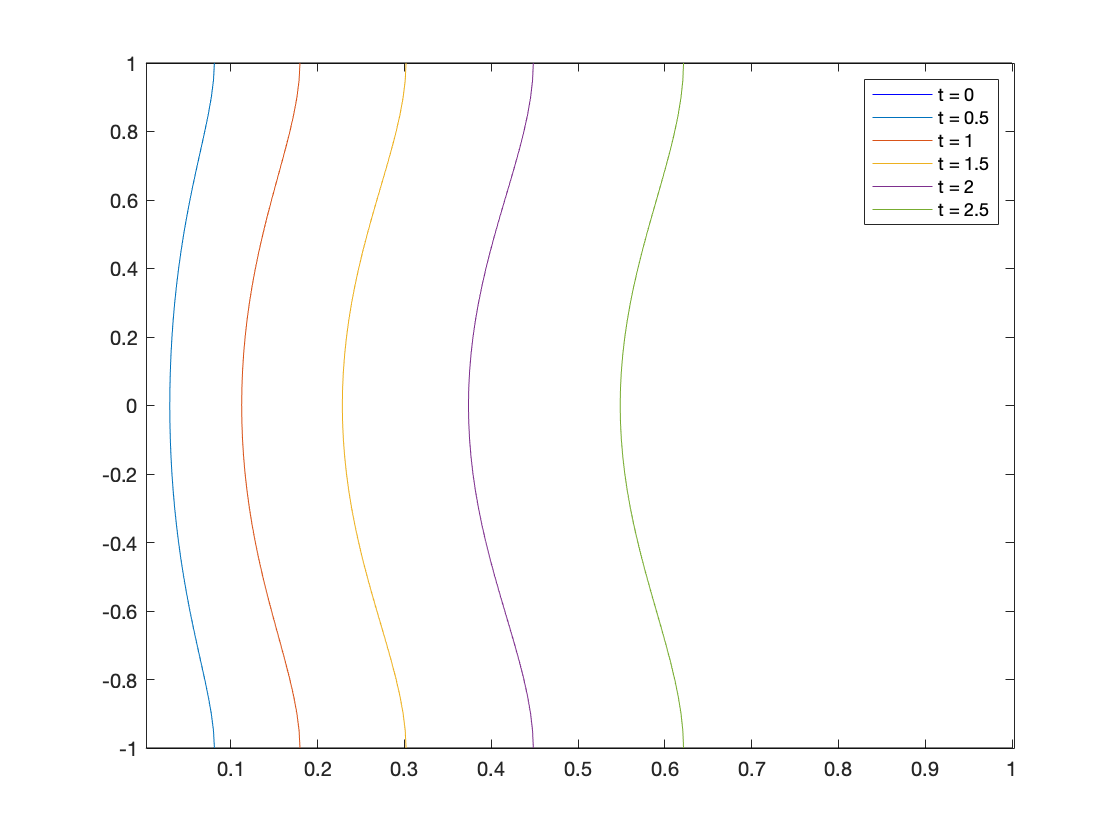}} 
~
\subfigure
{ \includegraphics[width=0.45\textwidth]{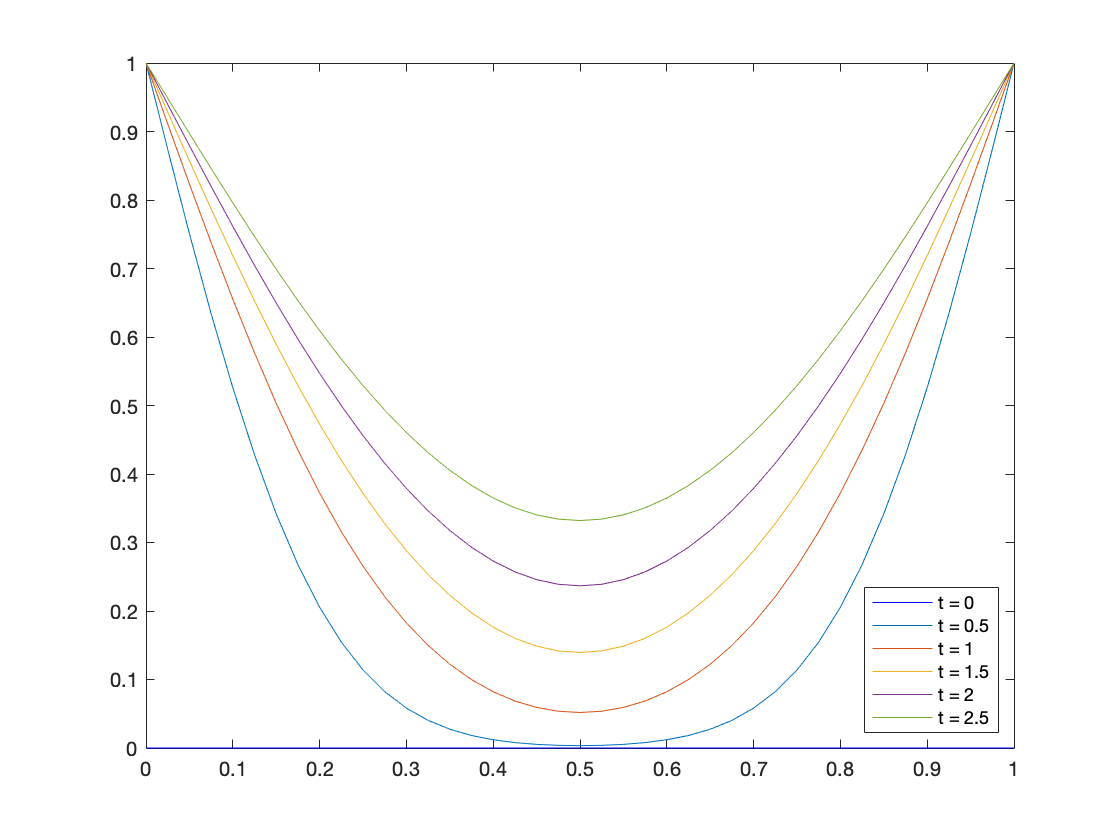}} 
\caption{DIGM simulation in a domain with straight boundaries.}
\label{fig:1}
\end{figure}

The set-up we consider here is similar to that considered in \cite{DES01}. Indeed the evolution law for the parametric system derived in \cite{DES01} can be obtained from \eqref{xeq} by setting $\alpha = 1$ and $F(\vp) = |\vp_0| - 1$, for some $\vp \in \mathbb{R}^2$, and considering a slightly different formulation of the reaction-diffusion equation \eqref{weq}. Setting $T = 2.5$ and
\begin{align*}
\vx^{\,0}(\rho) = (2 \rho - 1,0)^T, \quad w^0(\rho) = 0, \quad f(w) = w^2, \quad g(v,w) = v w, \qquad \rho \in [0,1]
\end{align*}
with the boundary data
\begin{align*}
F(\vp) = |\vp_0| - 1, \quad \gd = 1,
\end{align*}
yields the results displayed in Figure \ref{fig:1} in which a travelling wave solution is reached, the left hand plot shows the evolution of the interface at $t=0,0.5,1,1.5,2,2.5$, while the right hand plot shows the evolution of the solute, plotted against $\rho$, at the same times. These results are consistent with Figure 5--8 in \cite{DES01}.

{\bf Example 3}

\begin{figure}[ht!]
\centering
\subfigure
{ \includegraphics[width=0.45\textwidth]{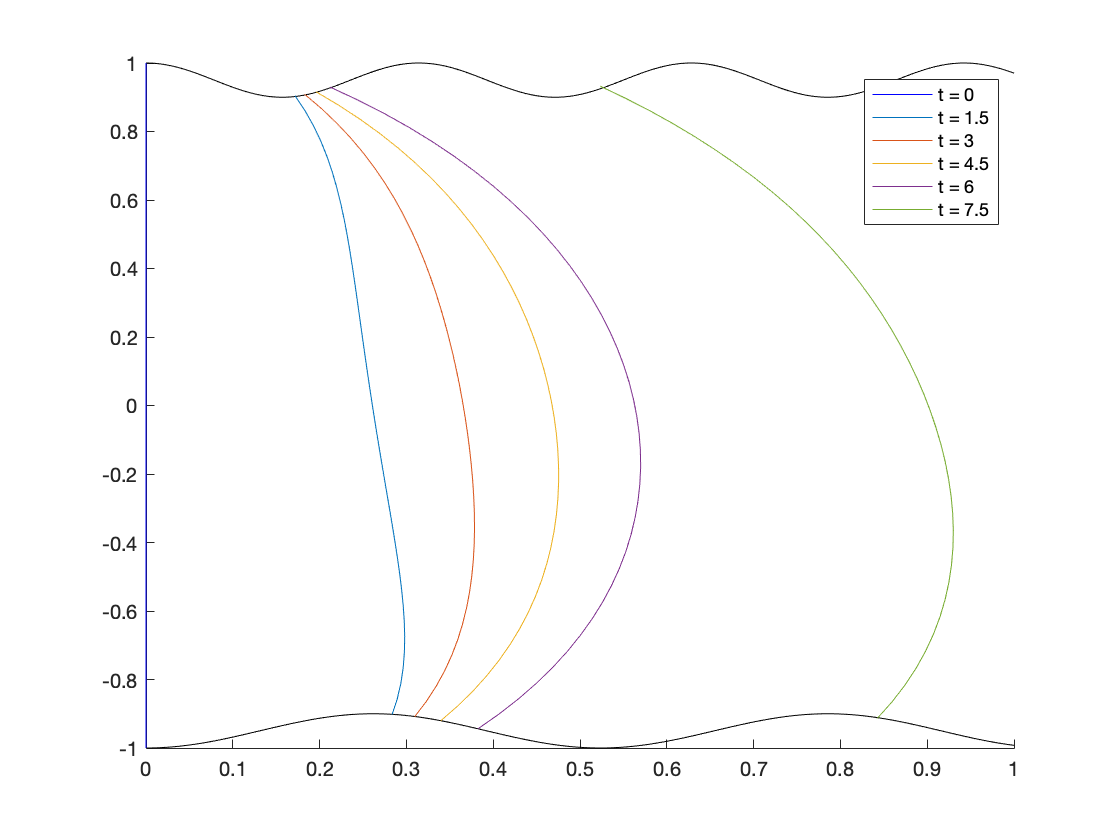}} 
\subfigure
{ \includegraphics[width=0.45\textwidth]{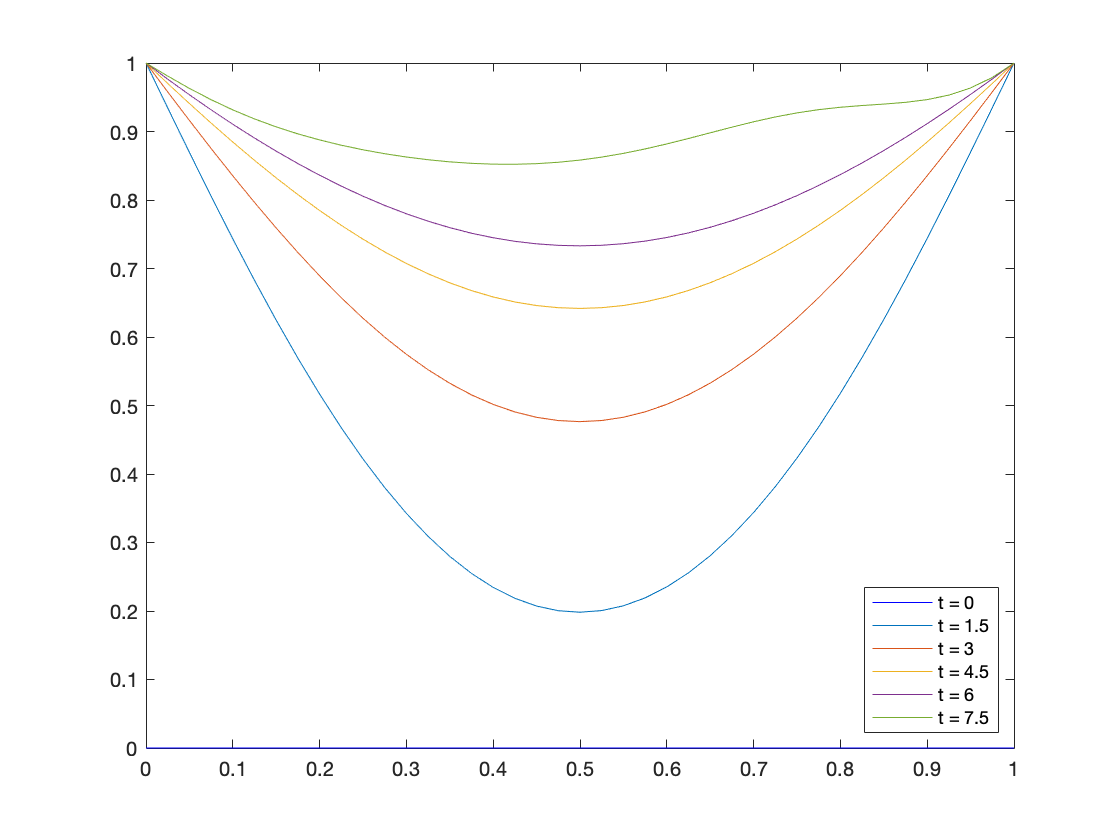}} 
\caption{DIGM simulation in a complex geometry.}
\label{fig:2}
\end{figure}

In this example we use the same data as Example 2 with the exception that we replace the simple straight-sided geometry of $\Omega$ that arose from setting $F(\vp) = |\vp_0| - 1$, with the more complex geometry
\begin{align*}
\Omega := \{ \vp \in \mathbb{R}^2 \text{ : } 0.05 \cos(20\vp_1) + 0.95 > \vp_0, \, -0.05 \cos(12\vp_1) - 0.5 < \vp_0 \}
\end{align*}
for which we note that \eqref{norm_F} does not hold.
The results are presented in Figure \ref{fig:2}, with the left hand plot displaying the evolution of the interface at $t=0,1.5,3,4.5,6,7.5$, together with the geometry $\Omega$ (black line) while the right hand plot shows the evolution of the solute, plotted against $\rho$, at the same times, with $T = 7.5$. From this figure we see that the complex nature of the domain destroys the travelling wave solution that was present in Example 2. 

{\bf Acknowledgements}
JVY gratefully acknowledges the support of the EPSRC grant 1805391. VS would like to thank the Isaac Newton Institute for Mathematical Sciences for support and hospitality during the programme {\it Geometry, compatibility and structure preservation in computational differential equations} when work on this paper was undertaken.
This work was supported by: EPSRC grant number EP/R014604/1.

\end{document}